%% file: paper.tex
\documentclass[final,3p,times]{elsarticle}

\usepackage{graphicx}

\usepackage{amssymb}
\usepackage{amsmath}
\usepackage{subfig}
\usepackage{graphicx}
\usepackage{url}
\usepackage{algorithm}
\usepackage{algpseudocode}
\usepackage{multirow}
\usepackage{wasysym}
\usepackage{marvosym}
\usepackage{color, xcolor, colortbl}
\usepackage{array, booktabs}
\usepackage{xspace}
\usepackage{enumerate}
\usepackage[shortlabels]{enumitem}

\usepackage{bm}
\input{macros}

\input{colours}

\newcolumntype{a}{>{\columncolor{light}}c}
\newcolumntype{b}{>{\columncolor{skyblue}}c}
\newcolumntype{d}{>{\columncolor{lightpurple}}c}

\makeatletter
\def\ps@pprintTitle{%
  \let\@oddhead\@empty
  \let\@evenhead\@empty
  \def\@oddfoot{\reset@font\hfil\thepage\hfil}
  \let\@evenfoot\@oddfoot
}
\makeatother
\begin{document}

\begin{frontmatter}



\title{Coarsening and parallelism with reduction multigrids for hyperbolic Boltzmann transport\tnoteref{crown}}
\author[AMCG]{S. Dargaville}
\ead{dargaville.steven@gmail.com}
\tnotetext[crown]{UK Ministry of Defence © Crown owned copyright 2024/AWE}
\author[AWE,AMCG]{R.P. Smedley-Stevenson}
\author[AMEC,AMCG]{P.N. Smith}
\author[AMCG]{C.C. Pain}
\address[AMCG]{Applied Modelling and Computation Group, Imperial College London, SW7 2AZ, UK}
\address[AWE]{AWE, Aldermaston, Reading, RG7 4PR, UK}
\address[AMEC]{ANSWERS Software Service, Jacobs, Kimmeridge House, Dorset Green Technology Park, Dorchester, DT2 8ZB, UK}
\begin{abstract}
Reduction multigrids have recently shown good performance in hyperbolic problems without the need for Gauss-Seidel smoothers. When applied to the hyperbolic limit of the Boltzmann Transport Equation (BTE), these methods result in very close to $\mathcal{O}(n)$ growth in work with problem size on unstructured grids. This scalability relies on the CF splitting producing an $\mat{A}_\textrm{ff}$ block that is easy to invert. We introduce a parallel two-pass CF splitting designed to give diagonally dominant $\mat{A}_\textrm{ff}$. The first pass computes a maximal independent set in the symmetrized strong connections. The second pass converts F-points to C-points based on the row-wise diagonal dominance of $\mat{A}_\textrm{ff}$. We find this two-pass CF splitting outperforms common CF splittings available in \textit{hypre}. 

Furthermore, parallelisation of reduction multigrids in hyperbolic problems is difficult as we require both long-range grid-transfer operators and slow coarsenings (with rates of $\sim$1/2 in both 2D and 3D). We find that good parallel performance in the setup and solve is dependent on several factors: repartitioning the coarse grids, reducing the number of active MPI ranks as we coarsen, truncating the multigrid hierarchy and applying a GMRES polynomial as a coarse-grid solver.

We compare the performance of two different reduction multigrids, AIRG (that we developed previously) and the \textit{hypre} implementation of $\ell$AIR. In the streaming limit with AIRG, we demonstrate 81\% weak scaling efficiency in the solve from 2 to 64 nodes (256 to 8196 cores) with only 8.8k unknowns per core, with solve times up to 5.9$\times$ smaller than the $\ell$AIR implementation in \textit{hypre}.
\end{abstract}

\end{frontmatter}
\section{Introduction}
\label{sec:Introduction}
The Boltzmann Transport Equation (BTE) considered in this work is the linear integro-differential equation that represents the transport and interaction of neutral particles (like neutrons and photons) within a media/material. Computing deterministic numerical solutions of the BTE is difficult, in part because in the limit of transparent media particles are transported across the domain; this is known as the ``streaming'' limit where the discretised system is hyperbolic. In the opposite limit, where interactions are very strong and the integral operator is dominant, the media scatters particles heavily in angle/energy (resulting in dense angle/energy blocks). Furthermore, in applications such as radiation transport in a nuclear reactor, the spatial dimension requires many elements (millions to billions) to build geometry-conforming meshes, with material properties varying in space/energy. Tens to hundreds of DOFs are typically applied in the angular dimension, with the energy dimension discretised with tens to hundreds of energy groups. Furthermore, the aforementioned dense angle/energy blocks in non-transparent problems result in the number of non-zeros in the discretised linear system growing non-linearly with angle refinement. 

These difficulties challenge both discretisation and solver technologies, with the two intrinsically linked. The space/angle/energy discretisation must remain stable in the presence of both heavy streaming and discontinuities caused by large variations in material properties, while also being amenable to highly-parallel, matrix-free iterative methods. Currently, there is only one combination of discretisation and solver technologies that practically satisfies both these constraints; in terms of discretisation in each dimension: Discontinuous Galerkin (DG) Finite Element Method (FEM) in space with upwinding, \sn (or low-order DG FEM) in angle and multigroup (P$^0$ DG FEM) in energy. For the iterative method, source iteration (a Richardson method) is typically used. This depends on the space/angle discretisation giving lower-triangular structure in the streaming/removal operator for each energy/angle block. This allows a single iteration of a sweep (a matrix-free Gauss-Seidel in space/angle) to exactly invert the block. Furthermore there are optimal parallel decompositions on structured spatial grids which allow sweeps to scale out to millions of cores on the largest HPC systems \cite{Baker1998, Adams2013}.

Unfortunately, sweeps do not scale well in parallel on general unstructured spatial grids. Previously \cite{Dargaville2024a, Dargaville2024} we presented an iterative method based on additive preconditioning that does not require sweeps. The key to this is the recent development of reduction multigrid methods which work well in hyperbolic problems \cite{Southworth2017, Manteuffel2019, Manteuffel2019a}. In \cite{Dargaville2024a} we developed a reduction multigrid known as approximate ideal restriction with GMRES polynomials (AIRG). We used a stabilised FEM method with a CG stencil to discretise in space, which reduces the memory consumption but does not give lower-triangular blocks in the streaming/removal operator. Thankfully AIRG does not depend on block or lower-triangular structure in our linear system, unlike the sweep-based method outlined above. This allows the scalable use of different space/angle discretisations for Boltzmann problems. 

We showed that using AIRG to solve the BTE in the streaming limit resulted in very close to $\mathcal{O}(n)$ growth in work with space/angle refinement on unstructured spatial grids. We also showed that the combination of additive preconditioning with one V-cycle of AIRG applied to the streaming/removal operator and one V-cycle of \textit{BoomerAMG} applied to a DSA type operator gave very close to $\mathcal{O}(n)$ growth in work in scattering problems. To be able to scalably solve the BTE in parallel, we therefore require that AIRG applied to either the streaming or streaming/removal operator scales. The streaming limit is the most difficult to solve and hence in this paper, we examine the performance of AIRG in parallel, applied to the streaming operator on unstructured grids. We also introduce a highly parallel CF splitting algorithm designed to produce good quality $\mat{A}_\textrm{ff}$ matrices on which reduction multigrids rely. 

Previous work on using reduction multigrids in parallel on the BTE is limited to that of \cite{Hanophy2020}, who used the AIR reduction multigrids known as local AIR ($\ell$AIR) and Neumann AIR (nAIR) implemented in \textit{hypre} applied to the streaming operator of a DG S$_n$ discretisation. When tested from $P=3$ to $P=4096$ cores on unstructured grids, the authors found the solve time weak scaled like $\log(P)^f$ with $f \approx 1.31$; this gives a weak scaling efficiency of 15\%.

In \cite{Dargaville2024a} we compared AIRG to $\ell$AIR and found that in serial, AIRG could produce close to scalable solutions with 3--5$\times$ less work than $\ell$AIR in the solve. The results in this work show this factor is retained in parallel. In general, reduction multigrids used on hyperbolic problems feature operators and grid hierarchies which differ from those found in elliptic problems, and hence we describe different strategies required for improving the parallel performance. 
\section{Discretisation}
\label{sec:Discretisation}
The first-order form of the time-independent BTE is given by
\begin{equation}
\bm{\Omega} \cdot \nabla \psif + \Sigma_\textrm{t} \psif = \int_{\bm{\Omega}'} \int_{E'} \Sigma_\textrm{s} (\bm{r}, \bm{\Omega}' \rightarrow \bm{\Omega}, E' \rightarrow E) \psifd \textrm{d}E' \textrm{d}\bm{\Omega}' + S_\textrm{e}(\bm{r}, \bm{\Omega}, E),
\label{eq:bte}
\end{equation}
where $\psif$ is the angular flux at spatial position $\bm{r}$, in direction $\bm{\Omega}$  and at energy $E$. The  material properties are characterised by the total and scatter macroscopic cross sections, $\Sigma_\textrm{t}$ and $\Sigma_\textrm{s}$, respectively, with $S_\textrm{e}$ the source term. In this work we consider solving mono-energetic problems in the streaming limit, so $\Sigma_\textrm{t}=\Sigma_\textrm{s}=0$. 

We give a summary of our sub-grid scale FEM discretisation below, for more details please see \cite{hughes_variational_1998, hughes_multiscale_2006, candy_subgrid_2008, buchan_inner-element_2010}. A sub-grid scale discretisation decomposes the solution as $\bm{\psi} = \bm{\varphi} + \bm{\theta}$, where $\bm{\varphi}$ is the solution on a ``coarse'' scale and $\bm{\theta}$ is a correction to the solution on a ``fine'' scale. In angle we discretise with a P$^0$ DG FEM with constant area azi/polar elements. We normalise the basis functions so that the angular mass matrix is the identity. The coarsest angular discretisation we use gives one constant basis function per octant, we denote this level 1, with subsequent levels given by splitting elements at the midpoint of the azi/cosine polar coordinates. We use a linear Continuous Galerkin (CG) formulation in space on the coarse scale and a linear Discontinuous Galerkin (DG) in space on the fine scale. The resulting linear system can be written as 
\begin{equation}
\begin{bmatrix}
\mat{A} & \mat{B} \\
\mat{C} & \mat{D} \\
\end{bmatrix}
\begin{bmatrix}
\bm{\Phi} \\
\bm{\Theta} \\
\end{bmatrix}
=
\begin{bmatrix}
\mat{S}_{\bm{\Phi}} \\
\mat{S}_{\bm{\Theta}} \\
\end{bmatrix},
\label{eq:SGS_full}
\end{equation}
where ${\bm{\Phi}}$ and ${\bm{\Theta}}$ are vectors containing the coarse and fine scale expansion coefficients, we denote the number of unknowns in ${\bm{\Phi}}$ as NCDOFs and in ${\bm{\Theta}}$ as NDDOFs, with NDOFs = NCDOFs + NDDOFs. The discretised discontinuous solution, $\bm{\Psi}$, is given by the addition of the coarse and fine scale solutions. The discretised source terms for both scales are $\mat{S}_{\bm{\Phi}}$ and $\mat{S}_{\bm{\Theta}}$. If we use a uniform angular discretisation, then the matrices $\mat{A}$ and $\mat{D}$ are the standard CG and DG FEM matrices that result from discretising \eref{eq:bte}, with $\mat{B}$ and $\mat{C}$ formed from integrating the continuous basis functions against the discontinuous (and vice versa).

We replace $\mat{D}$ with $\hat{\mat{D}}$, in which block-diagonal sparsity has been introduced by enforcing vacuum conditions on each DG element boundary (hence removing the DG jump terms) and only including self-scattering (in problems with scattering). With a uniform S$_n$ or P$^0$ DG FEM angular discretisation we have further block-diagonal sparsity, with $\hat{\mat{D}}$ having $3 \times 3$ (or $4 \times 4$ in 3D) blocks on the diagonal, given the lack of coupling between angles; we can invert this trivially. If instead we use angular adaptivity with P$^0$ DG elements then we invert $\hat{\mat{D}}$ without fill-in (i.e., using an ILU(0) factorisation), given angles can be coupled \cite{Dargaville2024}. We then have $\hat{\mat{D}}^{-1}$ and form a Schur complement for $\bm{\Phi}$. Our discontinuous discretised solution, $\bm{\Psi}$, is therefore computed by
\begin{enumerate}
\item Solving for the stabilised continuous solution $\bm{\Phi}$: $(\mat{A} - \mat{B} \hat{\mat{D}}^{-1} \mat{C}) {\bm{\Phi}} = \mat{S}_{\bm{\Phi}} - \mat{B} \hat{\mat{D}}^{-1} \mat{S}_{\bm{\Theta}}$.
\item Solving for the discontinuous correction $\bm{\Theta}$: $\bm{\Theta} = \hat{\mat{D}}^{-1} (\mat{S}_{\bm{\Theta}} - \mat{C} \bm{\Phi})$.
\item Combining the (interpolated) continuous solution and the discontinuous correction: $\bm{\Psi} = \bm{\Phi} + \bm{\Theta}$.
\end{enumerate}
Steps 2 and 3 are trivial. For Step 1 in problems without scattering, we can build an assembled version of $\mat{A} - \mat{B} \mat{D}^{-1} \mat{C}$ (i.e., the streaming/removal operator) with linear growth in the nnzs with space and angle refinement. We then use the algebraic multigrid (AMG) described below in \secref{sec:airg}, designed for hyperbolic problems to solve the linear system in Step 1. Importantly $\mat{A} - \mat{B} \mat{D}^{-1} \mat{C}$ is asymmetric, non-normal and when a S$_n$ or P$^0$ DG FEM uniform angular discretisation is used, the matrix is block diagonal (in angle). Each angle block has the stencil of a spatial CG discretisation and hence does not feature lower-triangular structure. 

For Step 1 in problems with scattering, $\mat{A} - \mat{B} \mat{D}^{-1} \mat{C}$ should never be explicitly assembled as the nnzs grows nonlinearly with angle refinement. Hence we would use the additively preconditioned iterative method described in \cite{Dargaville2024a} in scattering problems. This relies on assembling only the streaming/removal operator and using the aforementioned multigrid to inexactly invert this as an additive preconditioner.

Hence the main constraint on solving our discretised problem scalably, either in streaming or scattering problems, is inverting the streaming/removal operator. The next section discusses the multigrid that we use to do this. 
\section{Reduction multigrid}
\label{sec:reduction}
Reduction multigrids are formed with a coarse/fine (CF) splitting of the unknowns. If we consider a general linear system $\mat{A}\mat{x}=\mat{b}$ of size $n \times n$, with the number of fine and coarse points $n_\textrm{F} + n_\textrm{C} = n$, we can write the block splitting as 
\begin{equation}
\begin{bmatrix}
\mat{A}_\textrm{ff} & \mat{A}_\textrm{fc} \\
\mat{A}_\textrm{cf} & \mat{A}_\textrm{cc}
\end{bmatrix}
\begin{bmatrix}
\bm{x_\textrm{f}} \\
\bm{x_\textrm{c}}
\end{bmatrix} = 
\begin{bmatrix}
\bm{b_\textrm{f}} \\
\bm{b_\textrm{c}}
\end{bmatrix}.
\label{eq:air_two}
\end{equation}
Ideal prolongation and restriction operators are given by 
\begin{equation}
\mat{P} = 
\begin{bmatrix}
\mat{W} \\
\mat{I}
\end{bmatrix}, \quad
\mat{R} = 
\begin{bmatrix}
\mat{Z} & \mat{I}
\end{bmatrix}, \quad \textrm{where} \quad \mat{Z}=-\mat{A}_\textrm{cf} \mat{A}_\textrm{ff}^{-1}, \mat{W} = -\mat{A}_\textrm{ff}^{-1} \mat{A}_\textrm{fc},
\label{eq:prolong}
\end{equation}
Reduction multigrids are characterised by different approaches for computing approximations to $\mat{Z}$ (and/or $\mat{W}$) and the smoothers used, given $\mat{A}_\textrm{ff}^{-1}$ can be dense. The coarse-grid matrix is computed with $\mat{A}_\textrm{coarse}=\mat{R}\mat{A}\mat{P}$ and this process is applied recursively to generate a multigrid hierarchy.

The approach we describe below computes an approximation, $\hat{\mat{A}}_\textrm{ff}^{-1} \approx \mat{A}_\textrm{ff}^{-1}$ and hence, we form approximate operators, e.g., $\mat{Z}=-\mat{A}_\textrm{cf} \hat{\mat{A}}_\textrm{ff}^{-1}$. We also perform F-point smoothing with the same $\hat{\mat{A}}_\textrm{ff}^{-1}$. Hence the two (related) requirements for good performance in our reduction multigrid are: having the CF splitting produce a rapid coarsening that also gives a well-conditioned (and/or diagonally dominant) fine-fine block $\mat{A}_\textrm{ff}$ and forming a good sparse approximate inverse of $\mat{A}_\textrm{ff}^{-1}$. 
\subsection{AIRG}
\label{sec:airg}
Approximate ideal restriction with GMRES polynomials (AIRG) uses fixed-order GMRES polynomials to compute $\hat{\mat{A}}_\textrm{ff}^{-1}$ \cite{Dargaville2024a}. These polynomials give good approximations with symmetric or asymmetric $\mat{A}_\textrm{ff}$ while also being agnostic to unknown ordering, block and/or triangular structure. In the GMRES algorithm the solution at step $m$ can be written as $\mat{x}^m_\textrm{f} = q_{m-1}(\mat{A}_\textrm{ff}) \mat{b}_\textrm{f}$, where $q_{m-1}(\mat{A}_\textrm{ff})$ is a matrix polynomial of degree $m-1$ known as the GMRES polynomial and is given by
\begin{equation}
q_{m-1}(\mat{A}_\textrm{ff}) = \alpha_0 + \alpha_1 \mat{A}_\textrm{ff} + \alpha_2 \mat{A}_\textrm{ff}^2 + \ldots + \alpha_{m-1}\mat{A}_\textrm{ff}^{m-1},
\label{eq:gmres_poly}
\end{equation}
where $\alpha$ are the polynomial coefficients which minimise the 2-norm of the residual at step $m$ and can be output by a (slightly modified) GMRES method. Rather than use GMRES directly, in our multigrid setup we compute a stationary GMRES polynomial by generating a random rhs (normally distributed with zero mean) and then calculate our polynomial coefficients. We use this stationary polynomial as our approximate inverse, i.e., $\hat{\mat{A}}_\textrm{ff}^{-1} = q_{m-1}(\mat{A}_\textrm{ff})$. 

We found in \cite{Dargaville2024a} that low-order polynomials ($m \approx 3$) were sufficient to give good results with our multigrid, and hence we can use a communication-avoiding algorithm (by orthogonalising the power-basis in a single step) to compute the coefficients. We then explicitly built an an assembled matrix representation of our polynomial, with a matrix-matrix product used to compute our approximate ideal operators. We also used this assembled matrix as our F-point smoother during the solve; the stationary polynomials could alternatively be applied matrix-free if desired. In this work we only use assembled versions of our polynomials. 

In \cite{Dargaville2024a} we found this was a very effective approach for generating $\hat{\mat{A}}_\textrm{ff}^{-1}$. In an attempt to reduce the fill-in generated by the matrix powers, we also investigated constraining the sparsity of the matrix powers to have the same sparsity as $\mat{A}_\textrm{ff}$. This gives an approximate inverse as
\begin{equation}
\hat{\mat{A}}_\textrm{ff}^{-1} = \alpha_0 + \alpha_1 \mat{A}_\textrm{ff} + \alpha_2 \tilde{\mat{A}_\textrm{ff}^2} + \ldots + \alpha_{m-1} \tilde{\mat{A}_\textrm{ff}^{m-1}},
\label{eq:fixed_sparsity}
\end{equation}
where each subsequent matrix power (represented by tildes) has all entries outside of the sparsity of $\mat{A}_\textrm{ff}$ dropped. Fixing the sparsity in this manner gives a distance 2 approximate ideal restrictor. We found this fixed-sparsity polynomial very effective and cheap to assemble. In \cite{Dargaville2024a} we simply computed each of the matrix powers in \eref{eq:gmres_poly} using matrix-matrix products and then dropped entries outside the required sparsity. 

In this work we have instead written a fast method to compute the sum of fixed-sparsity powers of arbitrary order. We refer to enforcing the sparsity of $\mat{A}_\textrm{ff}$ as first order fixed sparsity, enforcing the sparsity of $\mat{A}_\textrm{ff}^2$ as second-order fixed sparsity, etc. Enforcing a fixed order of sparsity requires fewer local FLOPs and enables us to do far less communication in parallel. For example, for first order fixed sparsity, we gather the required entries in the off-process rows of $\mat{A}_\textrm{ff}$ once; the communication required for this is less than that needed to compute $\mat{A}_\textrm{ff}^2$. The numeric phase is able to compute the first order fixed-sparsity matrix powers with this data independently for each row.  Furthermore the local compute is very cache friendly as it operates on the same local and gathered non-local non-zero entries, regardless of the order of the polynomial. The only other communication required when computing the polynomial is that required to compute the polynomial coefficients, namely $m$ matvecs and a tall-skinny QR, which requires a single all-reduce, involving an $m \times m$ matrix in parallel. 

In summary, the setup of our AIRG multigrid requires on each level:
\begin{enumerate}
\item CF splitting
\item Extraction of required submatrices, e.g., $\mat{A}_\textrm{ff}$
\item Generation of a normally distributed rhs (with a Box-Muller transform), $\tilde{\mat{b}}_\textrm{f}$
\item $m$ matvecs to compute the power basis, $\mat{K}_{m+1} = [\tilde{\mat{b}}_\textrm{f}, \mat{A}_\textrm{ff} \tilde{\mat{b}}_\textrm{f}, \mat{A}_\textrm{ff}^2 \tilde{\mat{b}}_\textrm{f}, \ldots, \mat{A}_\textrm{ff}^{m} \tilde{\mat{b}}_\textrm{f} ]$
\item QR decomposition of $\mat{K}_{m+1}=\mat{Q}\mat{R}$ computed with a tall-skinny QR
\item Solution of the least-squares system $(\alpha_0, \ldots, \alpha_{m-1})^\textrm{T}= \textrm{argmin}_{y_m} ||\beta \mat{e}_1 - \tilde{\mat{R}} \mat{y}_m||_2$, where $\mat{e}_1$ is the first column of the $m+1$ identity, $\tilde{\mat{R}} = \mat{R}(:, {2}{:}{\textrm{end}})$ and $\beta = \mat{R}(1,1)$
\item Calculation of matrix powers of $\mat{A}_\textrm{ff}$ (with or without fixed sparsity applied)
\item Building an assembled version of $\hat{\mat{A}}_\textrm{ff}^{-1} = q_{m-1}(\mat{A}_\textrm{ff})$
\item Matrix-matrix product for the restrictor in \eref{eq:prolong}, $\mat{Z} = -\mat{A}_\textrm{cf} \hat{\mat{A}}_\textrm{ff}^{-1}$
\item Equivalent calculation for the prolongator or the generation of a classical one-point prolongator
\item Matrix-triple-product, $\mat{R} \mat{A} \mat{P}$ (or equivalent) to generate the coarse grid matrix
\end{enumerate} 
In the following section we describe the CF splitting algorithm that has been specifically tailored for parallel reduction multigrids. 
\section{CF splitting}
\label{sec:modified_pmis}
Classical multigrid coarsening algorithms (e.g., Ruge-Stuben) rely on heuristics to produce CF splittings that result in good grid transfer operators and small, sparse coarse grids (see \cite{Alber2007} for a review). Ruge-Stuben coarsening for example uses:
\begin{enumerate}[start=1,label={\bfseries H\arabic*:}]
\item For each node $j$ that strongly influences an F-point $i$, $j$ is either a C-point or strongly depends on a C-point $k$ that also strongly influences $i$.
\item No C-point strongly depends on another C-point.
\end{enumerate}
The heuristic \textbf{H1} ensures that the CF splitting produced and the Ruge-Stuben interpolation operators are compatible. For a reduction multigrid based on ideal operators, we have more freedom in the heuristics. The ideal grid-transfer operators, \eref{eq:prolong}, are automatically compatible with any given CF splitting; we exploit this property several times in this work. The only requirement for an effective reduction multigrid is that $\mat{A}_\textrm{ff}^{-1}$ is ``easy'' to sparsely approximate. For example, if the F-points are chosen to form an independent set in the adjacency graph of the matrix, then $\mat{A}_\textrm{ff}$ is diagonal and can be trivially inverted \cite{Saad1996, Saad1999, Saad2002}. 

For AIRG, we only require that $\mat{A}_\textrm{ff}^{-1}$ be easily approximated by a low-order GMRES polynomial. Tight convergence bounds for GMRES in non-normal systems are difficult to obtain (see \cite{Liesen2020} for example), but we found good performance from our reduction multigrid \cite{Dargaville2024a} when using the Falgout-CLJP algorithm in \textit{hypre} and low-order GMRES polynomials. Unfortunately with further experimentation we found the Falgout-CLJP algorithm is not performant in parallel on coarse grids in these problems. We also found that CF splitting algorithms with more parallelism, such as PMIS and HMIS, produce inferior coarsenings. The results from using these methods are given in \secref{sec:Results}. 

\cite{MacLachlan2007, MacLachlan2007a} produce CF splitting algorithms specifically designed to produce diagonally dominant $\mat{A}_\textrm{ff}$ (see also \cite{Saad1996}) that are well suited for use in reduction multigrids. They show that finding the largest F-point set such that every row of $\mat{A}_\textrm{ff}$ has a diagonal dominance ratio greater than some constant is NP-complete and hence they develop a greedy algorithm to compute an approximate solution. This greedy algorithm presented is serial as it processes nodes one at a time. \cite{Zaman2023, Zaman2022} use simulated annealing to approximately solve the same problem, though they find this approach is expensive. \cite{Reusken1999} compute CF splittings by finding a maximal independent set in the strong connections, though they use a breadth-first serial algorithm. 

We take inspiration from those works and seek a highly parallel CF splitting algorithm that produces a diagonally dominant $\mat{A}_\textrm{ff}$, whose inverse we assume can be well approximated by low-order matrix polynomials. We use a classical definition of strength of connection (SoC), where the set of nodes that node $i$ strongly depends on is denoted as $S_i$ and is given by
\begin{equation}
S_i = \{j:j \neq i, |a_{i,j}| \geq \alpha \max_{k \neq i} |a_{i,k}| \},
\label{eq:soc}
\end{equation}
where $\alpha$ is the strength threshold. Given we are solving asymmetric systems, we also need $S_i^\textrm{T}$ which are the strong influences of node $i$. We introduce two heuristics designed to encourage diagonally dominant $\mat{A}_\textrm{ff}$, namely
\begin{enumerate}[start=1,label={\bfseries H\arabic*R:}]
\item No F-point strongly depends on/influences another F-point.
\item Every row in $\mat{A}_\textrm{ff}$ should be strongly diagonally dominant.
\end{enumerate}
Enforcing \textbf{H1R} ensures there are no large off-diagonal entries in $\mat{A}_\textrm{ff}$, while also enforcing \textbf{H2R} ensures there are no rows with many small entries. \textbf{H1R} can be enforced by a Luby-type algorithm \cite{Luby1985}, where we seek a maximal independent set in the symmetrized strong connections. Specifically we use a parallel randomized greedy maximal independent set algorithm, which we denote as ``PMISR'' in Algorithm \ref{alg:pmisr}, to strongly enforce \textbf{H1R}.
\begin{algorithm}[ht]
\caption{Parallel Modified Independent Set for Reduction multigrids - PMISR}\label{alg:pmisr}
\begin{algorithmic}[1]
\State $C = \emptyset$
\ForAll{nodes $i$}
\State $w_i = S_i + S_i^\textrm{T} + \textrm{rand}(0,1)$
\EndFor
\State $F = \{i: w_i < 1\}$
\State $n_\textrm{loops} = 1$
\While{$n_\textrm{loops} > n_\textrm{loops}^\textrm{max}$ or there are nodes not in $F$ or $C$}
\State $D = \{i: w_i < w_j, \forall j \in S_i \cup S_i^\textrm{T}\}$
\ForAll{$j \in D$}
\State $F = F \cup j$
\ForAll{$k \in S_j \cup S_j^\textrm{T}$}
\State $C = C \cup k$
\EndFor
\EndFor
\State $n_\textrm{loops} = n_\textrm{loops} + 1$
\EndWhile
\State Any unassigned nodes are added to $C$
\end{algorithmic}
\end{algorithm}
This is very similar to the PMIS algorithm \cite{Sterck2006} with a symmetrized strength matrix and the definitions of the F and C points swapped. We should note, however, that they are not equivalent. This is because PMIS ensures each C-point has the largest measure amongst its strong dependencies/influences, whereas we want each F-point to have the smallest measure; in serial this is equivalent to modifying the order in which nodes are visited. Thankfully this is a trivial modification to existing PMIS implementations.  

After a first pass with Algorithm \ref{alg:pmisr}, we perform a second pass with Algorithm \ref{alg:cleanup} in order to weakly enforce \textbf{H2R}. We denote this ``cleanup'' phase as DDC and it converts any F-point with a row-wise diagonal dominance ratio, $\theta_i$, greater than some threshold, $\alpha_\textrm{diag}$, to a C-point. This is like a part of the greedy algorithm of \cite{MacLachlan2007, MacLachlan2007a} (or the preselection phase of \cite{Saad1996}), where instead we only perform a single iteration of this cleanup to remove the most problematic F-points. We also remove all the F-points greater than $\alpha_\textrm{diag}$ at once; this is more aggressive than optimal, as the removal of any one F-point can make other F-points more diagonally dominant in $\mat{A}_\textrm{ff}$ and hence they may otherwise be retained as F-points. Thankfully, as mentioned above, the automatic compatibility of a given CF splitting and the ideal operators means the only consequence of this is a slower coarsening. 
\begin{algorithm}[ht]
\caption{Diagonal dominance cleanup - DDC}\label{alg:cleanup}
\begin{algorithmic}[1]
\ForAll{$i \in F$}
\State $\theta_i = \frac{\sum_{j \in F, j \neq i} |a_{i,j}|}{|a_{i,i}|}$
\If{$\theta_i > \alpha_\textrm{diag}$ and $\theta_i \neq 0$}
\State $F = F \, \backslash \, i$
\State $C = C \cup i$
\EndIf
\EndFor
\end{algorithmic}
\end{algorithm}

\begin{figure}[th]
\centering
\subfloat[][Level 1]{\label{fig:diag_dom_compare}\includegraphics[width =0.4\textwidth]{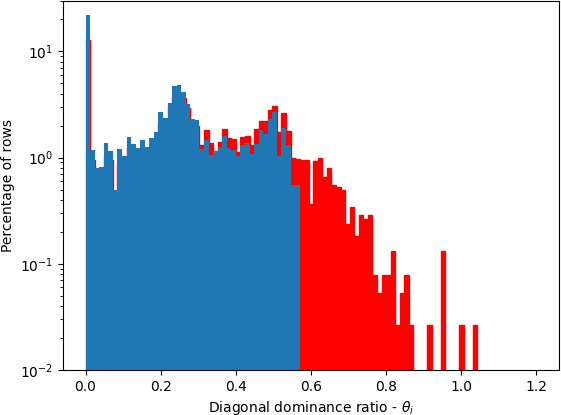}} \quad
\subfloat[][Level 2]{\label{fig:diag_dom_compare_level_2}\includegraphics[width =0.4\textwidth]{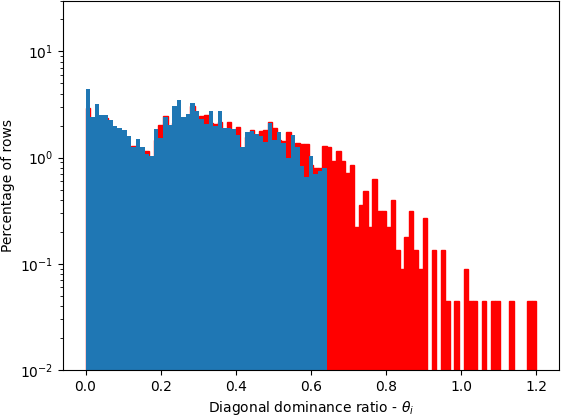}}
\caption{Row-wise diagonal dominance of $\mat{A}_\textrm{ff}$ on an unstructured mesh with 2313 nodes. \textcolor{red}{Red} shows after the first pass with PMSIR with strong tolerance 0.5, \textcolor{blue}{blue} shows after the second pass with DDC with fraction of 10\%.}
\label{fig:cleanup_hist}
\end{figure}

In summary, our CF splitting algorithm is one pass of PMISR followed by one pass of DDC. We have two parameters, $\alpha$ the strength parameter and $\alpha_\textrm{diag}$ the maximum row-wise diagonal dominance ratio. Choosing $\alpha_\textrm{diag}<1.0$ enforces strict diagonal dominance. We have found that rather than choosing $\alpha_\textrm{diag}$ directly, choosing $\alpha_\textrm{diag}$ such that a constant fraction of the least diagonally dominant F-points are converted to C-points performs well across a range of problems. The $\alpha_\textrm{diag}$ that gives this constant fraction can be calculated quickly by using the QuickSelect algorithm, or an approximate $\alpha_\textrm{diag}$ can be calculated by binning the diagonal dominance ratios and  choosing the closest bin boundary. We use this approach with 1000 bins. 
\section{Results}
\label{sec:Results}
To illustrate the performance of our CF splitting and AIRG, we show results from a streaming problem in 2D (zero scatter and absorption cross-sections). Our test problem is a 3$\times$3 box, with an isotropic source of strength 1 in a box of size 0.2$\times$0.2 in the centre of the domain. We apply vacuum conditions on the boundaries and discretise with unstructured triangles and ensure that the grids we use in any weak scaling studies are not semi-structured (e.g., we don't refine coarse grids by splitting elements). We use uniform level 1 refinement in angle, with 1 angle per octant (similar to S$_2$). Our matrix is therefore block diagonal, where each angle block is a stabilised CG discretisation of a time-independent advection equation (excluding boundary conditions), where the direction of the velocity is given by the angular direction. We do not exploit the block structure in this work; we apply the CF splitting and AIRG directly to the entire matrix. 

In the results below, we apply AIRG as a solver with an undamped Richardson with a relative tolerance of 1\xten{-10} and zero initial guess. We solve $(\mat{A} - \mat{B} \hat{\mat{D}}^{-1} \mat{C}) {\bm{\Phi}} = \mat{S}_{\bm{\Phi}} - \mat{B} \hat{\mat{D}}^{-1} \mat{S}_{\bm{\Theta}}$, which was denoted as Step 1 in \secref{sec:Discretisation}. We perform two iterations of up F-point smoothing, with no down-smooths. For the restrictor we use the approximate ideal restrictor shown in \eref{eq:prolong}, whereas for the prolongator we use a classical one-point operator \cite{Manteuffel2019}. 
\subsection{Serial comparison}
\label{sec:2d_stream}
In this section we show the results from a single unstructured spatial mesh to give some insight into the behaviour of our CF splitting algorithm, before moving into parallel. We use the same parameters for AIRG as in \cite{Dargaville2024a}; we should note, however, in this work that we use the SoC defined in \eref{eq:soc} with absolute values, so our results are not identical even when using the same CF splitting algorithm. To approximate $\mat{A}_\textrm{ff}^{-1}$ we use 3rd order GMRES polynomials ($m=4$) with fixed sparsity, with one iteration of a 3rd order GMRES polynomial with fixed sparsity used as the coarse grid solver. We terminate coarsening when we have fewer DOFs on the coarse grid than the GMRES polynomial order. We also drop small entries based on the relative row-wise infinity norm, using 0.0075 on the coarse grid matrices and 0.025 on the approximate ideal restrictor. 

\fref{fig:cleanup_hist} shows a histogram of the diagonal dominance ratios for each row of $\mat{A}_\textrm{ff}$ after our CF splitting algorithm, on the first two levels of the AIRG hierarchy when used on a problem with 2313 CG nodes (4464 elements) in space. We choose a strong tolerance of 0.5 to begin.  We can see that the first pass of PMISR results in a small number of rows that are not diagonally dominant, which is caused by the sum of small off-diagonal entries. We also see that there is a very sharp drop off in the diagonal dominance ratios, with the vast majority of rows strongly diagonally dominant. This motivates our approach of choosing a (small) fraction of F-points to convert to C-points in the DDC pass. \fref{fig:cleanup_hist} also shows the result of performing a second pass with DDC, with the worst 10\% of F-points changed to C-points. We can see that this reduces the worst dominance ratio from around 1.0 on the top grid and 1.2 on the second grid, to around 0.6 on both grids. 
\begin{figure}[th]
\centering
\subfloat[][Angle 1 - direction $(1,1)$]{\label{fig:angle_1}\includegraphics[width =0.3\textwidth]{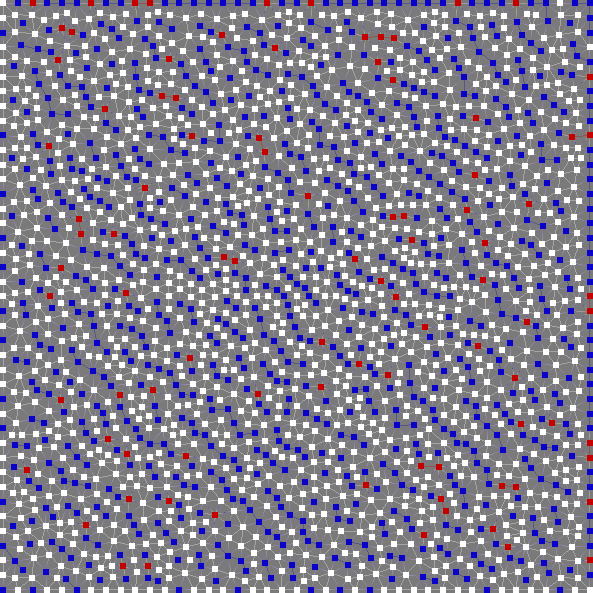}} \quad
\subfloat[][Angle 2 - direction $(-1,1)$]{\label{fig:angle_2}\includegraphics[width =0.3\textwidth]{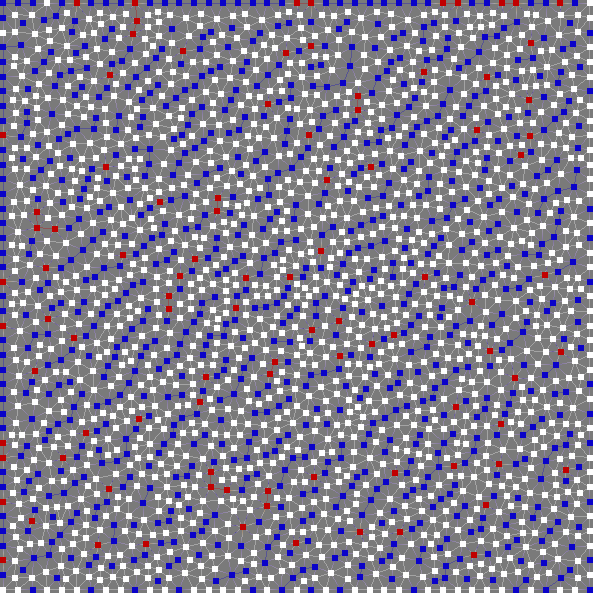}}
\caption{CF splitting produced by PMISR DDC on the top grid with strong tolerance 0.5 on an unstructured mesh with 2313 nodes. White squares are C-points, \textcolor{blue}{blue} squares are F points. The \textcolor{red}{red} squares are C-points that were converted from F-points by the second pass with DDC with fraction of 10\%.}
\label{fig:cf_splitting}
\end{figure}

\fref{fig:cf_splitting} shows the resulting C and F points on the top grid, for two different angles. We can see in both Figures \ref{fig:angle_1} and \ref{fig:angle_2} that a semi-coarsening is produced, with the F-points aligned (roughly) perpendicular to the angular direction (i.e., the ``flow''), as might be expected. Only weak off-diagonal entries should therefore be present in $\mat{A}_\textrm{ff}$. We can also see that the DDC second pass has converted F-points to C-points across the spatial mesh. We don't show the CF splitting for angles 3 and 4, though they are very similar to those shown in \fref{fig:cf_splitting} given our uniform angular discretisation and angular symmetry; we could exploit this and compute the CF splitting and/or symbolic matrix-matrix products for half the angles and then allow re-use, giving lower setup times (dependent on drop tolerances/padded sparsity). 

\subsubsection{Strong tolerances}
\label{sec:strong_tol}
We find PMISR DDC is robust across a wide range of strong tolerances. We show a comparison between three different strong tolerances in \tref{tab:strong_tol}, namely 0.0, 0.5 and 0.99. We measure the number of levels, the grid complexity (i.e., the number of unknowns across the hierarchy scaled by the number of unknowns on the top grid), operator complexity (i.e., the sum of nnzs in the coarse matrices across the hierarchy scaled by the top grid nnzs), iteration count and the total WUs (the number of FLOPs required to converge scaled by the FLOPs required to apply the top grid matrix). 

We also measure the ``storage complexity'', which reflects the relative storage required to perform a solve. With only up F-point smoothing we do not need to store the full matrix on each level once the next coarse matrix is produced in the setup, we only require $\mat{A}_\textrm{ff}$ and $\mat{A}_\textrm{fc}$. On the coarsest level we must store the approximate inverse and if we perform more than one iteration of our coarse solve we must also store the coarse matrix (to compute a residual). If $v_\textrm{coarse}$ is the number of iterations in the coarse solve and we have $l_{\textrm{max}}$ levels, then the storage complexity is
\begin{equation}
\textrm{Str. comp.} = \frac{\{\hat{\mat{A}}^{-1}\}^{l_{\textrm{max}}} + \min(v_\textrm{coarse}-1, 1)\{\mat{A}\}^{l_{\textrm{max}}} + \sum_{l=1}^{l=l_{\textrm{max}} - 1} \{\hat{\mat{A}}_\textrm{ff}^{-1}\}^l + \{\mat{A}_\textrm{ff}\}^l + \{\mat{A}_\textrm{fc}\}^l + \{\mat{R}\}^l + \{\mat{P}\}^l}{\{\mat{A}\}^1},
\label{eq:mat_store}
\end{equation}
where $\{.\}^l$ is the nnzs in an assembled matrix on level $l$.

\begin{figure}[th]
\centering
\subfloat[][Strong tolerance 0.99]{\label{fig:strong_099}\includegraphics[width =0.3\textwidth]{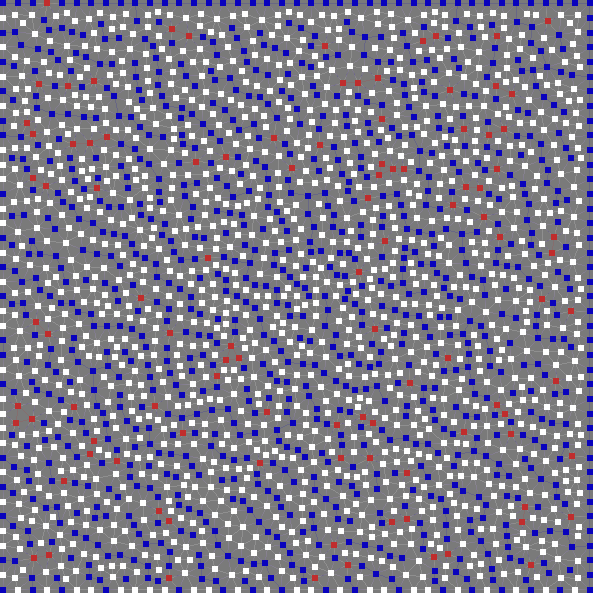}} \quad
\subfloat[][Strong tolerance 0.0]{\label{fig:perfect_is}\includegraphics[width =0.3\textwidth]{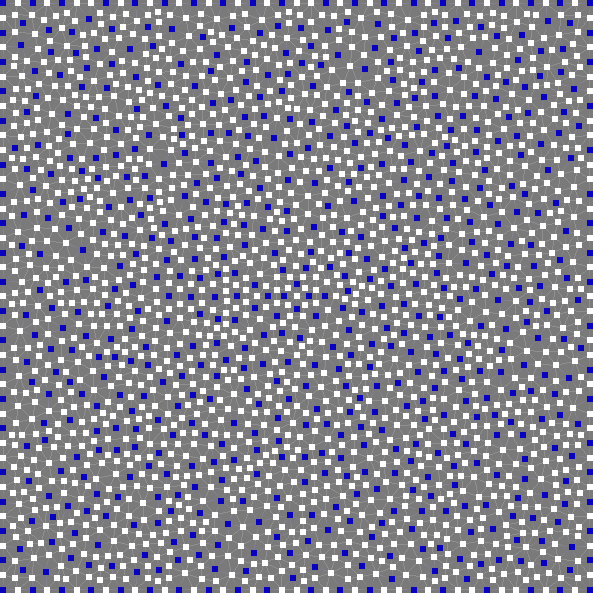}}
\caption{CF splitting produced for angle 1 (direction $(1,1)$) by PMISR DDC on the top grid of an unstructured mesh with 2313 nodes. White squares are C-points, \textcolor{blue}{blue} squares are F points. The \textcolor{red}{red} squares are C-points that were converted from F-points by the second pass with DDC with fraction of 10\%.}
\label{fig:strong_tol_fig}
\end{figure}

\tref{tab:strong_tol} shows that with a strong tolerance of 0.99 and one pass of PMISR, we achieve a grid complexity of 2.06. In advection problems (or a uniform angle discretisations of Boltzmann problems in the streaming limit), if we wish to avoid strong connections in $\mat{A}_\textrm{ff}$ we are fundamentally limited to coarsening along the characteristics. For example, in 2D on a structured quad mesh with grid aligned velocity, the ideal coarsening produces alternating lines of C and F points in the velocity direction, giving a grid complexity of exactly two. We cannot achieve such a perfect coarsening on an unstructured mesh, but we still see good convergence with 10 iterations; \fref{fig:strong_099} shows the CF splitting for angle 1 with a strong tolerance of 0.99 on our unstructured mesh. We see the coarsening is slower than in \fref{fig:angle_1}, but we still see a clear semi-coarsening. \tref{tab:strong_tol} shows that adding the second pass with DDC and a 10\% fraction slows the coarsening, with the grid complexity increasing from 2.06 to 2.3, but improves the convergence while also reducing the WUs from 59 to 53. 
\begin{table}[ht]
\centering
\begin{tabular}{c | c | c c c | c | c c }
\toprule
Strong tol. & CF splitting & Lvls & Grid comp. & Op. comp. & Str. comp. & its & WUs \\
\midrule
\multirow{2}{*}{0.99} & PMISR & 9 & 2.06 & 3.0 & 3.6 & 10 & 59 \\
& PMISR DDC & 10 & 2.3 & 3.4 & 3.7 & 9 & 53 \\
\midrule
\multirow{2}{*}{0.5} & PMISR & 12 & 2.5 & 3.5 & 3.5 & 8 & 43 \\
& PMISR DDC & 14 & 2.8 & 3.9 & 3.6 & 9 & 48  \\
\midrule
0.0 & PMISR & 39 & 6.0 & 9.2 & 3.9 & 7 & 34 \\
\bottomrule  
\end{tabular}
\caption{Performance of AIRG and quality measures given different strong tolerances on a streaming problem in 2D on an unstructured mesh with 2313 nodes.}
\label{tab:strong_tol}
\end{table}

We can see that decreasing the strong tolerance to 0.5, either with or without a second pass of DDC, produces a slower coarsening but also results in less work units compared to a tolerance of 0.99. In this case, the second pass of DDC has increased the number of WUs compared to a single pass of PMISR. As might be expected, DDC can ``cleanup'' a poor CF splitting and improve convergence, but if the CF splitting is good a second pass of DDC will slow down the coarsening without much benefit. In both the 0.99 and 0.5 case, the operator complexity and storage complexity are very similar as might be expected. 

Moving to the strong tolerance of 0.0, on the top grid this gives a maximal independent set in the connectivity of the spatial mesh for each angle, as the symmetrized strong connections with a strong tolerance of 0.0 is simply the adjacency graph of the matrix; see \fref{fig:perfect_is}. The benefit of producing an independent set in the adjacency graph of the matrix is that $\mat{A}_\textrm{ff}$ is diagonal. Only PMISR is shown in \tref{tab:strong_tol} with a strong tolerance of 0.0, as the DDC second pass ignores diagonal rows. This type of CF splitting and resulting diagonal $\mat{A}_\textrm{ff}$ have been investigated in several contexts, including cyclic reduction methods \cite{Heller1976, Golub1992, Gander1997}, reduction multigrid \cite{Ries1983} or early block LU/LDU preconditioners \cite{Saad1996}. The disadvantage is a substantially slower coarsening, with 39 levels instead of around 10 and a grid complexity of six.   

We can see, however, that although the operator complexity is 9.2, this is not a good surrogate for the memory required during the solve, as the storage complexity is 3.9. This is only slightly higher than with the two other tolerances. We see that an independent set coarsening actually gives a method with the best convergence, with only 7 iterations and 34 WUs. We should note we only perform one up F-point smooth with a (diagonal) fixed-sparsity GMRES polynomial as $\hat{\mat{A}}_\textrm{ff}^{-1}$, rather than just exactly inverting diagonal $\mat{A}_\textrm{ff}$ as we find better performance. These results suggest that using an independent set CF splitting may produce an efficient method, although the greater number of levels may increase the setup time, particularly in parallel; we examine this below. 

One of the benefits of an independent set CF splitting for Boltzmann transport problems is that with a uniform angle discretisation, the block structure in angle means that an independent set in the adjacency graph of the matrix is independent of angle (ignoring the boundary conditions). This allows even greater re-use of the CF splitting and symbolic matrix-matrix products than relying on angular symmetry as mentioned above (again dependent on drop tolerances/padded sparsity) and despite the greater number of levels this may lead to low setup times.

\subsubsection{Serial comparison of CF splitting algorithms}
\label{sec:serial_comparison}
We now compare common CF splitting algorithms in this serial problem before moving to parallel, as the small problem size means we can feasibly compute some (expensive) quality measures. \tref{tab:cond_aff} shows the results of using different CF splitting algorithms with AIRG, this time with a fixed strong tolerance of 0.5. For the Falgout-CLJP, PMIS, PMIS ``swap'' and HMIS algorithms we use the \textit{hypre} implementations. For both the Falgout-CLJP, PMIS and HMIS algorithms, we input the strength matrix as is, whereas for PMIS ``swap'', we input the symmetrized strength matrix to PMIS and the definition of F and C points are swapped on output. 

In this section we also show the maximum diagonal dominance ratio $\theta_i$, for each row of $\mat{A}_\textrm{ff}$ across all levels, along with the condition number, $\kappa$, of $\mat{A}_\textrm{ff}$ on each level. Furthermore, we show the maximum disk bound \cite{Liesen2020} across all levels, which gives the worst-case relative drop in the residual from applying our GMRES polynomial approximations of $\mat{A}_\textrm{ff}^{-1}$. These measures help quantify if the CF splittings are producing diagonally dominant $\mat{A}_\textrm{ff}$ that can be easily approximated by our GMRES polynomials. 

\begin{table}[ht]
\centering
\begin{tabular}{ c | c | c c | c | c | c c c c c c c c c c c c c c c}
\toprule
CF splitting & Grid comp. & its & WUs & $\max \frac{||\mat{r}^m||}{||\mat{r}^0||}$ & $\max \theta_i$ & \multicolumn{15}{c}{$\kappa(\mat{A}_\textrm{ff})$ on each level} \\
\midrule
Falgout-CLJP & 2.3 & 9 & 51 & 0.46 & 1.66 & 3.2 & 4.4 & 4.4 & 4.8 & 3.9 & 3.4 & 3.0 & 1.8 & & & & & & & \\
PMIS swap & 2.9 & 9 & 56 & 0.59 & 1.9 & 5.1 & 5.3 & 4.5 & 3.0 & 3.5 & 2.6 & 2.7 & 2.5 & 2.0 & 2.1 & 1.8 & 1.6 & & &\\
HMIS & 1.9 & 11 & 63 & 0.98 & 2.2 & 5.8 & 7.3 & 6.7 & 6.7 & 3.9 & 2.5 & & & & & & & & &\\
PMIS & 2.5 & 15 & 95 & 1.2 & 2.4 & 8.3 & 8.3 & 7.5 & 5.4 & 3.9 & 3.0 & 2.4 & 1.9 & 1.7 & 1.7 & 1.7 & 1.8 & 1.3 & &\\
\midrule
PMISR & 2.5 & 8 & 43 & 0.23 & 1.19 & 2.5 & 3.7 & 3.7 & 3.5 & 3.2 & 2.4 & 1.9 & 2.0 & 1.7 & 1.6 & 1.5 & 1.3 & & &  \\
PMISR DDC & 2.8 & 9 & 48 & 0.21 & 0.68 & 2.5 & 3.6 & 3.1 & 3.1 & 3.1 & 2.6 & 2.1 & 2.0 & 1.8 & 1.8 & 1.6 & 1.4 & 1.4 & 1.3 &  \\
\bottomrule  
\end{tabular}
\caption{Performance of AIRG and quality measures given different CF splittings with a strong threshold of 0.5 on a streaming problem in 2D on an unstructured mesh with 2313 nodes.}
\label{tab:cond_aff}
\end{table}

We can see in \tref{tab:cond_aff} that all the CF splitting algorithms perform well in these asymmetric problems, converging with less than 100 WUs and giving $\mat{A}_\textrm{ff}$ with $\mathcal{O}(1)$ condition numbers. Falgout-CLJP performs the best out of the \textit{hypre} implementations, with the lowest WUs. HMIS gives the most aggressive coarsening, followed by PMIS, Falgout-CLJP, with PMIS swap coarsening the slowest. We also see that the worst case convergence bounds seem to increase along with the worst diagonal dominance ratio; the closer to diagonal dominant the better the GMRES polynomial approximation. 

Furthermore \tref{tab:cond_aff} shows our PMISR CF splitting (without a second pass of DDC) performs the best out of all the CF splitting algorithms tested; the improved performance compared to PMIS swap shows that the two methods are not equivalent and hence it is worth ensuring each F-point has the smallest measure. The worst-case GMRES bound across the hierarchy with PMISR is roughly half that of the Falgout-CLJP method. We see that applying PMISR followed by a second pass of DDC with a 10\% fraction almost halves the worst diagonal dominance ratio, while only slightly decreasing the worst-case GMRES bound. In practice our GMRES polynomials are very strong approximations and we find they are less sensitive to the diagonal dominance ratio than say a Jacobi method. Overall we find either PMISR or PMISR DDC produce the best CF splittings, with the lowest condition numbers, diagonal dominance ratios and worst-case GMRES bounds. On more refined spatial meshes, typically we find that PMISR DDC performs the best. We now move to parallel results. 
\subsection{Parallel results}
\label{sec:parallel_results}
We now investigate the parallel performance of our CF splitting algorithm and AIRG multigrid. We show strong and weak scaling results in parallel on ARCHER2, a HPE Cray EX supercomputer with 2 $\times$ AMD EPYC 7742 64-core processors per node (128 cores per node). All timing results are taken from compiling our code, PETSc 3.14 and hypre with “-O3” optimisation flag. Our goal in this section is to achieve the smallest solve time, as we can amortize the setup cost over many solves in multi-group Boltzmann problems. Given this, we searched the parameter space to try and find the best values for drop tolerances, strength of connections, etc, although more optimal values may be possible. 

Both Algorithms \ref{alg:pmisr} \& \ref{alg:cleanup} are very parallel and suitable for GPUs. We can introduce two further modifications to improve the parallel performance. In PMISR we avoid the reduction at each iteration which computes the global number of unassigned nodes in Line 7 of Algorithm \ref{alg:pmisr}. Instead we allow the specification of a maximum number of iterations with $n_\textrm{loops}^\textrm{max}$, with any unassigned nodes made C-points; this means PMISR only requires nearest neighbour communication in parallel. \cite{Fischer2019} note that we would expect $\log(n)$ iterations in our parallel independent set algorithm, where in our case $n$ is the number of non-zeros in $S + S^\textrm{T}$. Given this log behaviour, the number of unassigned nodes decreases rapidly (and only a fraction of those would become F-points anyway). We find that even with hundreds of millions of unknowns, $n_\textrm{loops}^\textrm{max}=3$ is sufficient to see very little change in grid complexities. In DDC, we chose a constant fraction of the worst local F-points to convert to C-points, rather than the worst global F-points, in order to avoid an expensive global comparison. This makes the DDC pass dependent on the parallel decomposition, but we find this has a very limited effect and could be easily avoided by picking $\alpha_\textrm{diag}$ directly. 

In the results below, unless otherwise noted we use one pass of PMISR with a strong tolerance of 0.5 and $n_\textrm{loops}^\textrm{max}=3$ and one pass of DDC with local fraction 10\% as the CF splitting, along with a drop tolerance of 0.001 on the coarse matrix and 0.01 on the restrictor.
\subsubsection{GMRES polynomial order}
\label{sec:gmres_poly_order}
In our previous work \cite{Dargaville2024a}, we introduced a fixed-sparsity GMRES polynomial as an approximation to $\mat{A}_\textrm{ff}^{-1}$. We computed the fixed-sparsity matrix powers by calling the PETSc matrix-matrix product and dropping entries outside the sparsity of $\mat{A}_\textrm{ff}$. As noted in \secref{sec:airg} we have now written a fast method to compute fixed-sparsity matrix powers with arbitrary sparsity order. \tref{tab:poly_order_time} shows the performance of this method, with the time required on 8 nodes of ARCHER2 to compute the GMRES polynomial coefficients and assemble the matrix representing the approximate inverse (i.e., points 3--8 in \secref{sec:airg}) across all levels of the hierarchy (including the coarsest). This test was run with 4.5M elements and 2.3M CG nodes. 

We can see that as expected, the higher the order of the GMRES polynomial used, the more expensive it is to setup, but the better the convergence in the resulting multigrid. Moving from first order to third decreases the iteration count from 37 to 10, but the cost increases by approximately $19\times$. \tref{tab:poly_order_time} shows that the fixed-sparsity polynomials (i.e., sparsity $>$ order) continue to improve the convergence, while still being very cheap to compute. For example, with first order sparsity (i.e., the sparsity of $\mat{A}_\textrm{ff}$), increasing the polynomial order from 1 to 3 decreases the iteration count from 37 to 9, while only increasing the time from 0.037s to 0.092s. We see similar results for higher polynomial order. Choosing the sparsity order is a balancing act between the desired convergence, memory use and setup costs. Regardless of the sparsity order chosen, it is always beneficial to increase the polynomial order beyond this as we find there is not a significant cost to do so. Given this, in the remainder of this paper, unless otherwise noted we use 6th order GMRES polynomials ($m=7$) with first order fixed sparsity; hence we have a distance 2 approximate ideal restrictor.
\begin{table}[ht]
\centering
\begin{tabular}{c | a a a a | b b b b | d d d d}
\toprule
\rowcolor{white}
Order & Sparsity & its & WUs & Time (s) & Sparsity & its & WUs & Time (s) & Sparsity & its & WUs & Time (s) \\
\midrule
\rowcolor{white}
1 & 1 & 37 & 310 & 0.037 & & & & & & & & \\
2 & 1 & 17 & 137 & 0.084 & \cellcolor{white}2 & \cellcolor{white}15 & \cellcolor{white}140 & \cellcolor{white}0.27 & \cellcolor{white} &\cellcolor{white} &\cellcolor{white} &\cellcolor{white}\\
3 & 1 & 9 & 72 & 0.092 & 2 & 9 & 85 & 0.36 & \cellcolor{white}3 & \cellcolor{white}10 & \cellcolor{white}138 & \cellcolor{white}0.71\\
4 & 1 & 8 & 64 & 0.099 & 2 & 7 & 66 & 0.37 & 3 & 7 & 97 & 0.87\\
5 & 1 & 8 & 64 & 0.105 & 2 & 7 & 66 & 0.38 & 3 & 7 & 97 & 0.87\\
6 & 1 & 8 & 64 & 0.11 & 2 & 7 & 66 & 0.48 & 3 & 7 & 97 & 0.84\\
\bottomrule  
\end{tabular}
\caption{Time to compute and build assembled versions of GMRES polynomials across all levels on 8 nodes (1024 cores) of ARCHER2. The shaded rows represent fixed-sparsity GMRES polynomials (polynomial order $>$ sparsity order).}
\label{tab:poly_order_time}
\end{table}
\subsubsection{Optimisations for parallel multigrid}
\label{sec:parallel_multigrid}
Before we investigate the strong and weak scaling, we first discuss commonly used methods to improve the performance of multigrid methods in parallel. \cite{May2016} discussed three main strategies, namely: truncating the multigrid hierarchy, allowing a subset of MPI ranks to have zero unknowns and agglomerating unknowns onto a new MPI communicator with fewer MPI ranks. These can be used with topology aware communication strategies in the matrix-vector and/or matrix-matrix products \cite{Bhatele2008, Bienz2019, Bienz2020} along with repartitioned coarse grids \cite{Alef1991, Adams2004}.

Increasingly non-local operators can change the balance of work to communication on the coarse grids, while the natural difference in coarsening rates on different ranks with unstructured grids can alter the load balance. Repartitioning can manage both these trends, while also either retaining or decreasing the number of active MPI ranks. The repartitioning can be simple (e.g., odd ranks send their data to even ranks) or mesh partitioning software like ParMETIS \cite{karypis_metis-unstructured_1995} can be used to balance the work and communication over any number of active ranks. 

\begin{figure}[ht]
\centering
\includegraphics[width =0.45\textwidth]{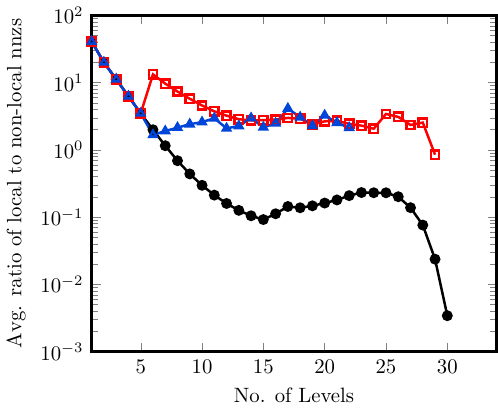}
\caption{Average ratio of local to non-local nnzs across active MPI ranks on each level on 32 nodes (4096 cores) of ARCHER2. Black is with no repartitioning, \textcolor{blue}{$\triangle$} is with simple repartitioning, \textcolor{red}{$\square$} is with ParMETIS repartitioning onto fewer ranks.}
\label{fig:nonlocal_ratio}
\end{figure}

\fref{fig:nonlocal_ratio} shows that without repartitioning, in a test with 18M elements and 9M CG nodes on 32 nodes of ARCHER2, our slow coarsening means we have around 30 levels in our hierarchy. We also see that the ratio of local to nonlocal work decreases as we coarsen, making our lower levels communication bound. We found this average ratio of local to nonlocal nnzs a good metric to trigger repartitioning; achieving good performance with this parameter is machine dependent. On ARCHER2 we use a value of 2. We chose to reduce the number of MPI ranks by half when we trigger repartitioning as this matches the expected coarsening rate of around 1/2 (in both 2D and 3D). 

\fref{fig:nonlocal_ratio} also shows a comparison between repartitioning with ParMETIS and a ``simple'' repartition where we simply halve the number of active MPI ranks (e.g., rank 1 sends all of its unknowns to rank 0, rank 3 to rank 2, etc.). We can see that both ways of repartitioning ensure a local to nonlocal ratio above two throughout the hierarchy. The ``simple'' repartitioning is triggered 10 times, giving lower grids on a single active MPI rank (hence the disappearance in \fref{fig:nonlocal_ratio} after level 23 as the work is entirely local). The ParMETIS repartitioning, however, is only triggered 5 times, as the minimisation of the edge cuts gives a significant improvement in the communication stencil required and hence the average local to nonlocal ratio; this is particularly visible after the first repartition on level 6. 

We could allow the repartitioning with ParMETIS to occur independently from the decrease in the number of active MPI ranks; for example, GAMG in PETSc allows the repartioning of every coarse grid, with the number of active MPI ranks only decreasing where necessary to ensure each local coarse grid has a set number of unknowns (the default is 50). We found combining these reduced the number of times we repartition, decreasing setup times without any large effect on solve times; on 32 nodes for example, we found repartioning every coarse level with ParMETIS and only decreasing the number of active ranks when the local to nonlocal nnzs ratio hit 2 actually increased our solve time by around 2-5\%. This is because of the extra communication introduced into the grid-transfer operators between levels when repartitioning has occurred. 

We use an ``interlaced'' repartitioning as we find this performs best, e.g., if our MPI rank numbering increases with the node number and the number of active MPI ranks is halved, each node has half the active number of ranks, rather than half the number of nodes being full. This also helps decrease the peak memory use on any given node. Furthermore, an interlaced repartition leaves idle cores on each node which allows the use of shared memory parallelism with OpenMP on lower levels. We experimented with this and found we could achieve good speedups in many aspects of the setup, particularly the assembly of the fixed-sparsity GMRES polynomials (Point 8 in \secref{sec:airg}). The matrix-matrix products used to form the restrictors and coarse matrices, however, did not see much speedup, as the use of drop tolerances in our hierarchy mean our matrices rarely become dense enough in the streaming limit; we found we required ${>}30$ nnzs per row to see good speedups. We did not explore the use of OpenMP parallelism during the solve on lower levels after repartioning, as it was difficult to implement the required non-busy MPI waits within the matrix-vector products used. We believe this would be worth exploring in the future, as the repartioning means the lower grid matrix-vector products potentially have enough local rows to see good speedups. We should note that none of the results in this paper use OpenMP in either the setup or the solve. 
\subsubsection{Timing results for parallel optimisations}
\label{sec:timing_parallel_optim}
Below we show results timing the solve and setup from the test problem in the previous section on 32 nodes. We examine a combination of the strategies listed above, namely no repartitioning, ``simple'' and ParMETIS repartitioning and the use of the node-aware matrix-vector and matrix-matrix products from the Raptor  library \cite{BiOl2017, Bienz2020}. We compare the use of two different sparse matrix-vector (SpMV) routines used on every level as part of the smoothing, grid-transfer operators and coarse grid solve; namely the PETSc and the Raptor SpMVs. Similarly we tested using the PETSc, \textit{hypre} and Raptor sparse matrix-matrix product (SpGEMMs) when computing the approximate ideal restrictor and the coarse grid matrices on every level. We found the \textit{hypre} SpGEMM was more expensive than both the PETSc and Raptor SpGEMMs and hence we don't show those results. The coarse grid matrix is computed with two SpGEMMs, i.e., $\mat{R} (\mat{A} \mat{P})$. We also enabled the ``scalable'' PETSc SpGEMM to ensure consistent performance. With Raptor, we found the best performance from setting the environmental variable \textsc{PPN} equal to 16; this is the number of cores per NUMA region on ARCHER2.
\begin{figure}[thp]
\centering
\subfloat[][No repartition, PETSc SpMV \& SpGEMM. \\Total setup: 42s, total solve: 0.26s]{\label{fig:setup_times}\includegraphics[width =0.4\textwidth]{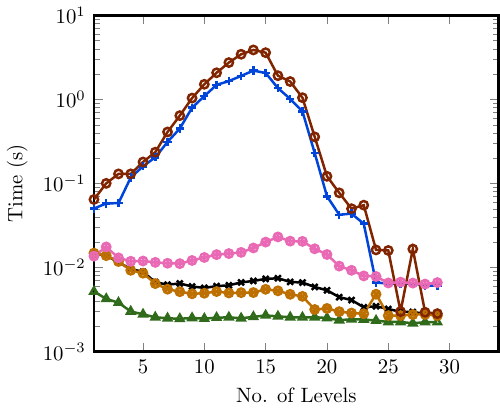}} \quad
\subfloat[][No repartition, Raptor SpMV \& SpGEMM. \\Total setup: 7.8s, total solve: 0.20s]{\label{fig:setup_times_raptor}\includegraphics[width =0.4\textwidth]{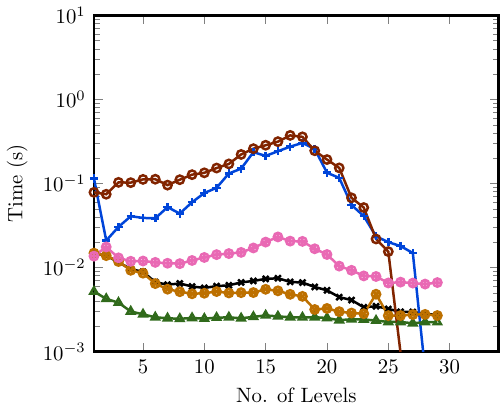}} \quad
\subfloat[][Simple repartition onto fewer ranks, PETSc SpMV \& SpGEMM. \\Total setup: 15.1s, total solve: 0.33s]{\label{fig:setup_times_simple}\includegraphics[width =0.4\textwidth]{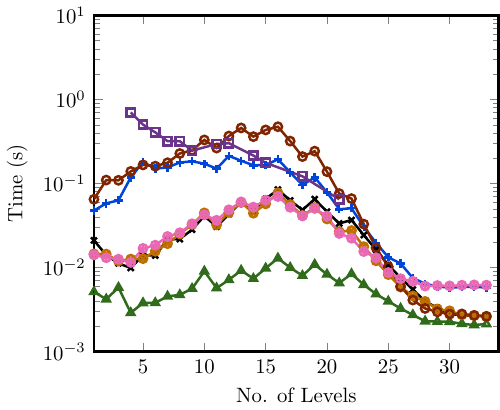}} \quad
\subfloat[][ParMETIS Repartition onto fewer ranks. Raptor SpMV \& SpGEMM. \\Total setup: 12.1s, total solve: 0.2s]{\label{fig:setup_times_parmetis_raptor}\includegraphics[width =0.4\textwidth]{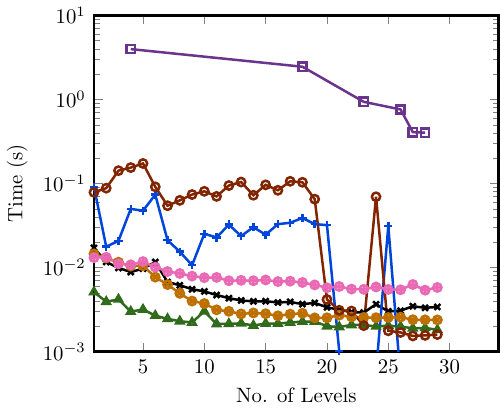}} \\
\subfloat[][ParMETIS Repartition onto fewer ranks, PETSc SpMV \& SpGEMM. \\Total setup: 13.0s, total solve: 0.086s]{\label{fig:setup_times_parmetis}\includegraphics[width =0.4\textwidth]{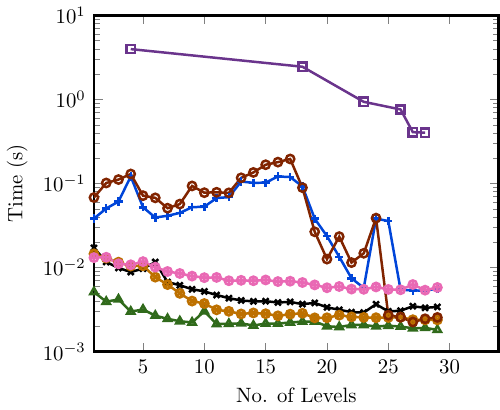}} 
\caption{Time taken for each component of the multigrid setup on each level on 32 nodes (4096 cores) of ARCHER2. The \textcolor{black}{$\times$} is the CF splitting, the \textcolor{foliagegreen}{$\triangle$} is the prolongator, the \textcolor{deludedorange}{$\otimes$} is the GMRES polynomial, the \textcolor{matlabblue}{$+$} is the SpGEMM for the restrictor, the \textcolor{fireenginered}{o} is the SpGEMM for the coarse grid, the \textcolor{darksalmon}{$\oplus$} is the matrix extract and the \textcolor{gaylordpurple}{$\square$} is the repartitioning.}
\label{fig:setup_times_whole}
\end{figure}

\fref{fig:setup_times_whole} shows timings for each component of our multigrid setup, as listed in \secref{sec:airg}. In this figure and the remainder of the paper, we split the cost of computing our approximate ideal restrictor into two parts: the cost of computing and assembling the GMRES polynomial (points 3--8 in \secref{sec:airg}) and the cost of the SpGEMM (point 9 in \secref{sec:airg}). 

We can see in \fref{fig:setup_times} that with no repartitioning, the SpGEMMs associated with our restrictor and the calculation of the coarse grid matrix are very expensive, particularly in the middle of the hierarchy where we are communication bound; the total setup time is 42 seconds, with the solve time 0.26 seconds. With no repartitioning, \fref{fig:setup_times_raptor} shows, however, that using the Raptor SpGEMM is significantly faster in the middle of the hierarchy. This gives a total setup time of 7.8s and a smaller solve time with the Raptor SpMVs of 0.2 seconds. 

The simple repartitioning shown in \fref{fig:setup_times_simple} also decreases the cost of the SpGEMMs significantly, while increasing the cost of computing both the GMRES polynomial and the CF splitting. This is because by increasing the local to non-local ratio, the simple repartitioning also increases the number of local rows and the cost of both the GMRES polynomial and CF splitting are closely tied to the local nnzs. We still find overall a significant decrease in setup time, at 15.1 seconds. The solve time, however, increases to 0.33 seconds. This is due to several effects: the increased number of local rows on active ranks after repartitioning, the increased communication required by the grid-transfer operators each time we hit a level that has been repartitioned (which happens 10 times with the simple repartitioning), along with the lack of optimisation of that communication. 

The combination of repartitioning with ParMETIS and using Raptor gives a reduced setup cost, at 12.1s; this is largely due to (roughly) halving the time required by the SpGEMM for the restrictor, as the time to compute the coarse matrix was approximately the same. The biggest expense in the setup is now the repartioning with ParMETIS, taking a total of 8.9 seconds compared with 3.6 seconds with the simple method; the solve time is at 0.2s. \fref{fig:setup_times_parmetis} shows that using ParMETIS repartioning but without Raptor, we have a slight increase in setup time at 13.0 seconds, but importantly the solve time has decreased to 0.086 seconds, which is the lowest value and 3$\times$ less than with no repartioning.
\begin{figure}[ht]
\centering
\includegraphics[width =0.45\textwidth]{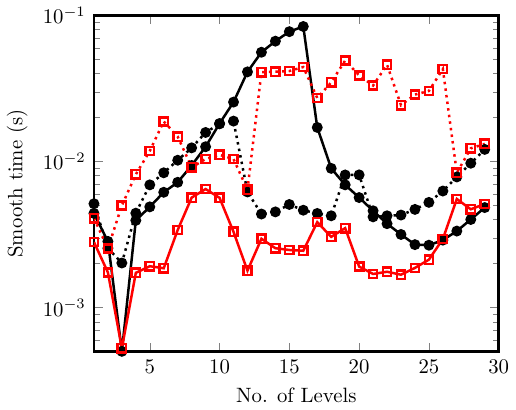}
\caption{Time spent in the F-point smooth on each level on 32 nodes (4096 cores) of ARCHER2. $\CIRCLE$ are with no repartitioning, \textcolor{red}{$\square$} are with ParMETIS repartitioning onto fewer ranks. The solid lines are the PETSc SpMV , the dotted are the Raptor SpMV.}
\label{fig:smooth_time}
\end{figure}

In order to further examine the impact of repartitioning and node-aware communication on the solve time, \fref{fig:smooth_time} shows the time spent in F-point smoothing in the solve across the hierarchy; the time spent in SpMVs for the grid-transfer operators behaves similarly. We can see that without repartionining and using the PETSc SpMV, we have a significant increase in smooth time in the middle of the hierarchy, much like with the PETSc SpGEMMs discussed above. We can also see that the Raptor SpMV reduces this cost in the middle of the hierarchy; it is this effect that lead to the reduced solve time reported in \fref{fig:setup_times_raptor}. Once repartitioning with ParMETIS has been enabled, however, we find that the Raptor SpMV can be more expensive on lower levels than the PETSc SpMV. We also compared the NAP-2 and NAP-3 communication options in Raptor but could not find settings that reduced the Raptor SpMV time below that of the PETSc SpMV when repartitioning was enabled. This is perhaps to be expected; node-aware communication is most beneficial when the communication volume is poorly optimised. 

Hence we conclude that repartitioning with ParMETIS is crucial to the performance of both the solve and setup. In the setup of this test problem, we found that the Raptor SpGEMMs was sometimes the cheapest. If we chose not to repartition, using the Raptor SpGEMM significantly reduced the setup time, giving the lowest value. If we repartition with ParMETIS, however, the Raptor SpGEMM was cheaper for some levels/operators and not others. We also found this in other problems when strong/weak scaling; the Raptor SpGEMM ranged from twice to half the cost of the PETSc SpGEMM when ParMETIS repartitioning was used. In the solve, we find it is only worth using the Raptor SpMV if we have not repartitioned with ParMETIS; the lowest solve time comes from repartitioning with ParMETIS and using the PETSc SpMVs. Given this, we would only recommend using the node-aware SpMV and SpGEMM from the Raptor if the coarse grids have not been repartitioned onto fewer MPI ranks with an optimising library such as ParMETIS. 

We can further improve the setup and solve time by truncating our multigrid hierarchy, but this requires a strong coarse solver. We found using an exact LU factorisation was far too expensive and resulted in poor scaling, in both the setup and solve phase. Instead we use our GMRES polynomials as an iterative coarse grid solver, truncating the hierarchy when the coarse grid has approximately 300k unknowns. We use three iterations of a 12th order GMRES polynomial with first order fixed sparsity. We also investigated using a 12th order GMRES polynomial applied matrix free (i.e., without first order fixed sparsity), but we found it faster to apply several iterations of the fixed-sparsity polynomial (at the cost of storing an assembled approximate inverse with fixed sparsity). 

The GMRES polynomial is very cheap to setup when used a coarse grid solver. For the example in \fref{fig:setup_times_parmetis}, truncating to 14 levels results in a cumulative GMRES polynomial setup time of 0.14s on levels 1--13 with the coarse level setup on level 14 costing just 0.03s. In the remainder of this paper, unless otherwise noted we truncate our hierarchy and use a GMRES polynomial as the coarse grid solver, use the PETSc SpMV and SpGEMM and also repartition our coarse grids with ParMETIS and half the number of active MPI ranks when the local to non-local nnzs ratio hits 2.
\subsubsection{Strong scaling}
\label{sec:strong_scaling}
\begin{table}[ht]
\centering
\begin{tabular}{c c c | c | c c c | c }
\toprule
Nodes & Cores & CDOFs/core & Solve time (s) & WU ratio & C-Grind time (ns) & D-Grind time (ns) & Setup time (s) \\
\midrule
2 & 256 & 35k & 0.152 & 1.17 & 540 & 90 & 3.8 \\
4 & 512 & 18k & 0.089 & 1.44 & 627 & 104 & 3.4 \\
8 & 1024 & 8.8k & 0.057 & 1.49 & 813 & 135 & 4.5 \\
20 & 2560 & 3.5k & 0.037 & 2.13 & 1322 & 220 & 5.4 \\
30 & 3840 & 2.3k & 0.032 & 2.81 & 1707 & 284 & 7.5 \\
\bottomrule  
\end{tabular}
\caption{Strong scaling of AIRG with ParMETIS repartitioning onto fewer ranks and a truncated hierarchy (13 levels) on ARCHER2, on a test problem with 4.5M elements and 2.3M CG nodes and 4 angles. The solve takes 8 iterations with 64 WUs.}
\label{tab:strong_scaling}
\end{table}
We now examine strong scaling results with the PMISR DDC CF splitting and AIRG. \tref{tab:strong_scaling} shows the results of a strong scaling experiment, run with 4.5M elements and 2.3M CG nodes. We find strong scaling efficiencies of 86\% on 4 nodes, 46\% on 20 and 40\% on 30 in the solve when compared to the 2 node case. The minimum wall time in the solve of 0.032s occurs with 30 nodes, down at 2.3k CDOFs per core, although this comes at the cost of increasing the setup time. As is the case in many algebraic multigrids, the setup time strong scales poorly and actually increases with high numbers of nodes. This is due to the increasing cost of the matrix-matrix products on lower grids, similar to that seen in \secref{sec:parallel_multigrid}. We should note that the repartioning with ParMETIS is vital to achieving the strong scaling in \tref{tab:strong_scaling}; without it on 30 nodes for example, the solve and setup take roughly 3$\times$ longer, at 0.099s and 24s, respectively.

The ``WU ratio'' in \tref{tab:strong_scaling} is the solve time divided by the WUs, scaled by the time to compute a top grid matvec. A value of 1 would indicate that the WUs in the solve perfectly predicts the observed solve time and hence there is no communication bottleneck on the lower grids. We can see that the WU ratio increases as we strong scale, up to 2.81 on 30 nodes. Similar to the setup, this is due to the increasing communication bottleneck on the lower grids. Again the repartitioning with ParMETIS helps keep this ratio low; without repartitioning on 30 nodes we have a WU ratio of 10.2.

\tref{tab:strong_scaling} also includes two different `grind times''. This is the total solve time, divided by either the NCDOFs or the NDDOFs per core and divided by the iteration count. This gives the cost per V-cycle scaled by the number of continuous or discontinuous DOFs per core. In \cite{Dargaville2024a} we show good results from using a single V-cycle on the streaming/removal operator as a preconditioner in problems with scattering. Hence we show the grind times as they quantify the cost of a single V-cycle. We can see from the C-Grind time that our solve is most efficient with more unknowns per core. At approximately 35k CDOFs/core we require 540 ns per CDOF/core/V-cycle. If we compare this with the grind time of a structured grid DG-sweep code (which is the time required per DDOF/core/sweep), this seems expensive; high-performance structured grid DG-sweeps can achieve times as low as 30-100 ns. We should note, however, that the C-Grind time is that required to perform a V-cycle on a stabilised \textit{continuous} discretisation, namely Step 1 in \secref{sec:Discretisation}. Our discretisation is discontinuous, though it relies on computing the solution of a stabilised continuous problem first, with the further work required to then compute our discontinuous solution small and trivially parallelisable (Step 2 and 3 in \secref{sec:Discretisation}). 

Hence we also also compute the D-Grind time, which shows the cost of solving our stabilised continuous problem but scaled by the number of discontinuous DOFs. This helps us make comparisons with DG-sweep codes; we achieve 90 ns per DDOF/core/V-cycle on an unstructured grid on ARCHER2. Of course a DG-sweep exactly inverts the streaming operator (ignoring the impact of cycles on unstructured grids), whereas the D-Grind time is that for a single V-cycle, which only reduces the residual by roughly an order of magnitude. As we showed in \cite{Dargaville2024a}, however, a single V-cycle is sufficient when the streaming/removal operator is used as a preconditioner. Using a single V-cycle as a preconditioner on the angular flux requires more memory in the solve when compared to the source iteration used in typical DG-sweep codes. When used with a stabilised continuous discretisation, however, the smaller number of DOFs helps reduce the required memory in the solve and results in grind times that are competitive with structured grid DG sweeps, even on unstructured grids. 
\subsubsection{Weak scaling}
\label{sec:weak_scaling}
The previous section established that around 35k CDOFs per core with AIRG gives good performance. The mesh infrastructure we have built our implementation on has a strict limit on the number of elements we can use; in parallel this limits the weak scaling results we can show. Given this, we decided to decrease the number of DOFs per core in our weak scaling experiments. This is in an attempt to increase the communication bottlenecks which would happen naturally at higher core counts, despite mesh limitations preventing us from running with a higher number of cores. Hence in the experiments below we show results from only using around 8.8k CDOFs per core. Similarly, we only use 4 angles in our angular discretisation (uniform level 1 refinement) in order to reduce the amount of local work, again exposing communication bottlenecks. 

We begin by showing weak scaling results without ParMETIS repartitioning onto fewer ranks or truncation of the multigrid hierarchy to give a baseline. \tref{tab:weak_scaling} shows that despite the iteration count only increasing from 8 to 9 (and the WUs growing from 62 to 71), the solve time goes up by a factor of over two with the setup time growing considerably from 2 to 64 nodes. \secref{sec:parallel_multigrid} showed this was due to the communication bottlenecks on the lower levels. If we use ParMETIS repartitioning onto fewer ranks, both the solve and setup time decrease considerably; on 64 nodes the solve time is 3$\times$ less. Truncating the hierarchy and relying on the strength of our GMRES polynomials as a coarse grid solver further reduces the solve and setup time as shown in \tref{tab:weak_scaling_repartion_truncate}, but without affecting the iteration count. We see a weak scaling efficiency of 76\% in the solve from 2 nodes to 64, with a plateau in solve time on 32 and 64 nodes. 

If we further truncate the hierarchy and increase the number of iterations used in the coarse solver, \tref{tab:weak_scaling_repartion_truncate_heavy} shows we achieve the same iteration count. This gives growth in the time per iteration of only 9.6\% from 2 nodes to 64. Overall we have 81\% weak scaling efficiency in the solve, even with very few DOFs per core. This comes at the cost of a more expensive solve, with an increase from 0.075s to 0.088s on 64 nodes (and an increase in the WUs required). In general this is further confirmation that although ParMETIS repartitioning onto fewer ranks is significantly improving the performance of our method, there is still some growth in solve time driven by communication limitations on the lower levels that must be managed. 

\begin{table}[ht]
\centering
\begin{tabular}{c c c | c c c | c | c c | c c }
\toprule
Nodes & Eles & CG nodes & Lvls & Grid comp. & Op. comp. & Str. comp. & its & WUs & Solve (s) & Setup (s) \\
\midrule
2 & 1.1M & 567k & 24 & 2.9 & 6.9 & 5.4 & 8 & 62 & 0.11 & 4.6 \\
8 & 4.5M & 2.2M & 27 & 2.9 & 7.1 & 5.5 & 8 & 63 & 0.15 & 16.5 \\
32 & 18M & 9M & 30 & 2.9 & 7.2 & 5.5 & 9 & 71 & 0.2 & 45.8 \\
64 & 36M & 18M & 31 & 2.9 & 7.2 & 5.5 & 9 & 72 & 0.23 & 64.8 \\
\midrule
Weak scaling efficiency & \multicolumn{7}{c}{} & & 47\% & 7.1\% \\
\bottomrule  
\end{tabular}
\caption{Weak scaling of AIRG on ARCHER2 without repartitioning or truncation (with 128 cores per node).}
\label{tab:weak_scaling}
\end{table}

\begin{table}[ht]
\centering
\begin{tabular}{c c c | c c c | c | c c | c c }
\toprule
Nodes & Eles & CG nodes & Lvls & Grid comp. & Op. comp. & Str. comp. & its & WUs & Solve (s) & Setup (s) \\
\midrule
2 & 1.1M & 567k & 11 & 2.8 & 6.7 & 5.4 & 8 & 64 & 0.052 & 2.2 \\
8 & 4.5M & 2.2M & 13 & 2.8 & 7.0 & 5.5 & 8 & 64 & 0.057 & 4.5 \\
32 & 18M & 9M & 14 & 2.8 & 7.1 & 5.5 & 9 & 72 & 0.067 & 6.8 \\
64 & 36M & 18M & 15 & 2.9 & 7.1 & 5.5 & 9 & 72 & 0.068 & 11.6 \\
\midrule
Weak scaling efficiency & \multicolumn{7}{c}{} & & 76\% & 19\% \\
\bottomrule  
\end{tabular}
\caption{Weak scaling of AIRG with ParMETIS repartitioning onto fewer ranks and a truncated hierarchy on ARCHER2 (with 128 cores per node).}
\label{tab:weak_scaling_repartion_truncate}
\end{table}

\begin{table}[ht]
\centering
\begin{tabular}{c c c | c c c | c | c c | c c }
\toprule
Nodes & Eles & CG nodes & Lvls & Grid comp. & Op. comp. & Str. comp. & its & WUs & Solve (s) & Setup (s) \\
\midrule
2 & 1.1M & 567k & 6 & 2.5 & 5.5 & 5.3 & 8 & 232 & 0.072 & 0.7 \\
8 & 4.5M & 2.2M & 7 & 2.6 & 6.0 & 5.4 & 8 & 192 & 0.083 & 4.2 \\
32 & 18M & 9M & 8 & 2.7 & 6.3 & 5.5 & 9 & 179 & 0.095 & 5.8 \\
64 & 36M & 18M & 9 & 2.7 & 6.6 & 5.5 & 9 & 151 & 0.088 & 7.9 \\
\midrule
Weak scaling efficiency & \multicolumn{7}{c}{} & & 81\% & 8.5\% \\
\bottomrule  
\end{tabular}
\caption{Weak scaling of AIRG with ParMETIS repartitioning onto fewer ranks and a heavily truncated hierarchy (using 20 iterations of the coarse solver instead of 3) on ARCHER2 (with 128 cores per node).}
\label{tab:weak_scaling_repartion_truncate_heavy}
\end{table}
\subsubsection{Weak scaling with strong tolerance 0.0}
\label{sec:weak_scaling_strong_zero}
\begin{table}[ht]
\centering
\begin{tabular}{c c c | c c c | c | c c | c c }
\toprule
Nodes & Eles & CG nodes & Lvls & Grid comp. & Op. comp. & Str. comp. & its & WUs & Solve (s) & Setup (s) \\
\midrule
2 & 1.1M & 567k & 49 & 8.5 & 23.1 & 5.9 & 6 & 42 & 0.067 & 3.0 \\
8 & 4.5M & 2.2M & 57 & 8.8 & 24.4 & 6.0 & 7 & 49 & 0.08 & 5.8 \\
32 & 18M & 9M & 78 & 9.01 & 25.5 & 6.1 & 7 & 49 & 0.11 & 9.8 \\
64 & 36M & 18M & 90 & 9.08 & 25.8 & 6.1 & 8 & 56 & 0.11 & 15.1 \\
\midrule
Weak scaling efficiency & \multicolumn{7}{c}{} & & 59\% & 19\% \\
\bottomrule 
\end{tabular}
\caption{Weak scaling of AIRG with strong tolerance 0.0 with ParMETIS repartitioning onto fewer ranks and a truncated hierarchy on ARCHER2 (with 128 cores per node).}
\label{tab:weak_scaling_repartion_truncate_zero}
\end{table}

In \secref{sec:strong_tol} we discussed the use of a strong tolerance of 0.0, giving a maximal independent set in the adjacency graph of the matrix on each level. \tref{tab:weak_scaling_repartion_truncate_zero} shows that in parallel this still results in a method with better convergence and less WUs than seen with a strong tolerance of 0.5 as in \secref{sec:weak_scaling}. The cost of this is a very slow coarsening, with between 49--90 levels and a grid complexity of between 8.5--9.08. The operator complexity is very high at around 25, but we can see that the storage complexity is comparable to that with higher strong tolerances, requiring around 6 copies of the top grid matrix. We can also see that the resulting solve and setup times are around 50\% higher than with a strong tolerance of 0.5 in \tref{tab:weak_scaling_repartion_truncate}, but they are still reasonably low. In a typical elliptic multigrid, a method with such slow coarsening and hence high grid and operator complexity would not be performant. 

It is the combination of only F-point smoothing and the ParMETIS repartitioning that is responsible for such good performance despite the slow coarsening. The small number of F-points means only F-point smoothing is cheap and doesn't require an expensive C-point residual or C-point smooth, giving a low number of WUs. Similarly the SpGEMMs required to compute the approximate ideal restrictor are very rectangular, with the SpGEMMs required to compute the coarse matrix benefiting from large identity blocks in $\mat{R}$ and $\mat{P}$ given the large number of C-points. Using a local to non-local nnzs ratio to trigger the ParMETIS repartitioning onto fewer ranks along with a slow coarsening also ensures that repartitioning only occurs a few times throughout the hierarchy, similar to with a higher strong tolerance. The repartitioning also helps minimise the communication bottlenecks in the solve on the (many) lower levels. We know that the strong tolerance of 0.5 is, however, less communication bound in the solve than with 0.0, as we see a reduction in WUs from 63 to 49 on 8 nodes for example, but an increase in solve time from 0.057s to 0.08s. The resulting WU ratio shows this, increasing from 1.59 with a strong tolerance of 0.5 (see \tref{tab:strong_scaling}) to 2.78 with a strong tolerance of 0.0 on 8 nodes. Without the ParMETIS repartitioning onto fewer ranks, however, both the solve and setup time in this case shows much poorer weak scaling. 

Using a strong tolerance of 0.0 and hence a diagonal $\mat{A}_\textrm{ff}$ may give a good method for streaming Boltzmann problems as the independent set means the coarsening is independent of angle. As mentioned in \secref{sec:strong_tol}, this may allow considerable re-use in the setup, with only a reasonable performance penalty in the solve in parallel as shown in \tref{tab:weak_scaling_repartion_truncate_zero}. We leave exploration of re-use and the impact of different strong tolerances to future work. 
\subsubsection{Parallel comparison of CF splitting algorithms}
\label{sec:parallel_comparison}
Similar to \secref{sec:serial_comparison} we show a comparison between different CF splitting algorithms, but now in parallel. We compare weak scaling results for the \textit{hypre} implementations of Falgout-CLJP, PMIS ``swap'' as defined in \secref{sec:serial_comparison} and HMIS against PMISR and PMISR DDC with 10\% local fraction. We don't show results from PMIS as it performed very poorly. In all these tests we use a strong tolerance of 0.5.

We can see in \tref{tab:parallel_cf_compare} that the lowest solve and setup time across all CF splitting algorithms tested occur with PMISR DDC. We can also see that the time required to compute the PMISR and PMISR DDC CF splitting is substantially smaller than the other algorithms. We also ran a comparison where we do not limit the maximum number of steps in our PMISR algorithm, i.e., we don't set $n_\textrm{loops}^\textrm{max}$. On 32 nodes with PMISR for example, the grid complexity with $n_\textrm{loops}^\textrm{max}=3$ is 2.588, whereas with no maximum the grid complexity is 2.583. This small difference makes no impact on the iteration count (or solve time), but setting $n_\textrm{loops}^\textrm{max}=3$ cuts the time spent in the CF splitting by 25\% on 32 nodes. 

Of the \textit{hypre} CF splitting algorithms, \tref{tab:parallel_cf_compare} shows that the Falgout-CLJP produces the best solve time, only 5\% slower than that of PMISR DDC. The time required by PMISR DDC to compute the CF splitting, however, is much less; on 32 nodes it is 23$\times$ faster. This confirms that our PMISR DDC CF splitting algorithm is very performant in parallel and an excellent choice for use with reduction multigrids which require diagonally dominant $\mat{A}_\textrm{ff}$. 
\begin{table}[ht]
\centering
\begin{tabular}{ c | c c c | c c c | c c c }
\toprule
CF splitting & \multicolumn{3}{c}{Solve (s)} & \multicolumn{3}{c}{CF splitting (s)} & \multicolumn{3}{c}{Setup (s)} \\
 & 2 nodes & 8 nodes & 32 nodes & 2 nodes & 8 nodes & 32 nodes & 2 nodes & 8 nodes & 32 nodes \\
\midrule
Falgout-CLJP & 0.07 & 0.088 & 0.089 & 0.56 & 2.5 & 6.7 & 3.3 & 8.5 & 18.1 \\
PMIS swap & 0.077 & 0.097 & 0.13 & 1.2 & 4.4 & 12.2 & 4.7 & 11.8 & 27.9 \\
HMIS & 0.11 & 0.13 & 0.15 & 3.3 & 0.7 & 1.5 & 4.4 & 8.3 & 14.4 \\
\midrule
PMISR & 0.0068 & 0.083 & 0.093 & \cellcolor{light} 0.14 & 1.0 & 0.97 & 3.6 & 8.0 & 14.2 \\
\rowcolor{light}
\cellcolor{white}PMISR DDC & 0.066 & 0.08 & 0.086 & \cellcolor{white} 0.16 & 0.27 & 0.29 & 2.8 & 7.6 & 13.0 \\
\bottomrule  
\end{tabular}
\caption{Weak scaling of AIRG with different CF splitting algorithms with strong tolerance 0.5, no truncation and ParMETIS repartitioning onto fewer ranks on ARCHER2 (with 128 cores per node). The minimum times are shaded.}
\label{tab:parallel_cf_compare}
\end{table}
\subsubsection{Weak scaling comparison with \textit{hypre} $\ell$AIR}
\label{sec:weak_scaling_compare}
To finish our weak scaling experiments, we compare our results with the \textit{hypre} implementation of distance 2 $\ell$AIR. We tried to use parameters that gave the lowest solve time and best weak scaling in order to give a fair comparison with AIRG. We use 3 FCF-Jacobi up smooths, 5 damped Jacobi iterations as a coarse solver and a strong $\mat{R}$ tolerance of 0.0. The default strong $\mat{R}$ tolerance is larger and we have found values greater than $0.0$ results in less WUs in the solve and smaller setup times for a given mesh. Using a larger tolerance, however, significantly degrades convergence and scaling (see also \cite{Dargaville2024a}). We used a Falgout-CLJP CF splitting with a strong tolerance of 0.2 as \secref{sec:parallel_comparison} (and \secref{sec:serial_comparison}) showed this gave the lowest solve time out of the algorithms available in \textit{hypre}. \secref{sec:parallel_comparison} showed the Falgout-CLJP algorithm was very expensive in parallel, so we should note the setup times shown below for $\ell$AIR include this. In future work we would like to examine the role of the $\ell$AIR approximate ideal restrictor separate from the implementation and CF splitting algorithms available in \textit{hypre}. 

\textit{hypre} does not provide the ability to repartition coarse grids and hence the only technique we have to reduce the communication bottlenecks on lower grids in these hyperbolic problems is truncation. We attempted to use $\ell$AIR without truncating the hierarchy in order to compare against \tref{tab:weak_scaling}, but we found the setup was too expensive. Hence to begin, we truncated the hierarchy and used 12, 15, 18 and 21 levels, respectively. This resulted in 16 iterations in the solve and weak growth in the WUs with 316, 326, 331 and 333 WUs. 

\begin{figure}[thp]
\centering
\subfloat[][Solve]{\label{fig:weak_scaling_solve}\includegraphics[width =0.4\textwidth]{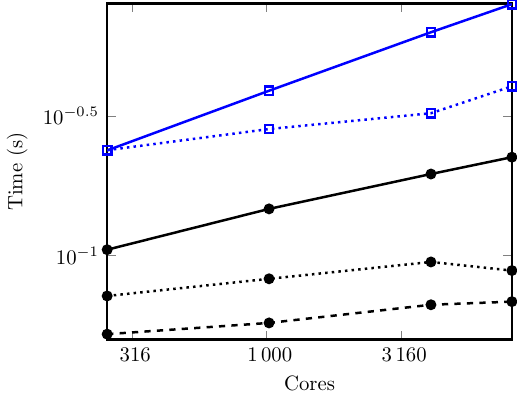}} \quad
\subfloat[][Setup]{\label{fig:weak_scaling_setup}\includegraphics[width =0.4\textwidth]{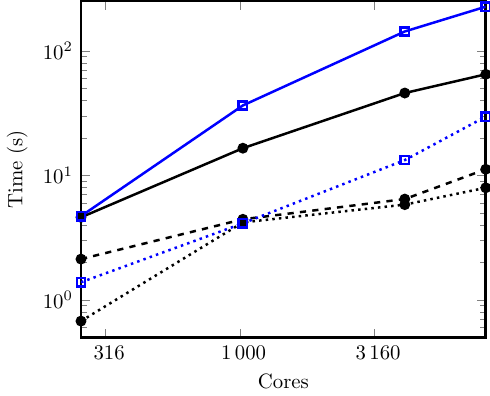}} 
\caption{Comparison of weak scaling results for AIRG and \textcolor{blue}{$\ell$AIR} on ARCHER2, from 2 nodes (256 cores) to 64 nodes (8196 cores) with 8.8k CDOFs/core. The solid \textcolor{blue}{$\square$} is $\ell$AIR with a truncated hierarchy, the dotted \textcolor{blue}{$\square$} is $\ell$AIR with a hevaily truncated hierarchy. The solid $\CIRCLE$ is AIRG with no truncation, the dashed $\CIRCLE$ is AIRG with ParMETIS repartitioning onto fewer ranks and a truncated hierarchy, the dotted $\CIRCLE$ is AIRG with ParMETIS repartitioning onto fewer ranks and a heavily truncated hierarchy.}
\label{fig:weak_scaling}
\end{figure}

We should note it is difficult to compare the WUs between AIRG and the \textit{hypre} implementation of $\ell$AIR as the cycle complexity output by \textit{hypre} doesn't include all the work associated with smoothing. In \cite{Dargaville2024a} we redefined the cycle complexity calculation in \textit{hypre} to enable comparisons, but it is more difficult to deploy this modification on ARCHER2. As such, we only note the WUs for $\ell$AIR to make clear that the heavy growth in the solve time shown in \fref{fig:weak_scaling_solve} is caused by communication bottlenecks, rather than growth in work. We can also see in \fref{fig:weak_scaling_solve} that $\ell$AIR with truncation is more expensive in the solve than AIRG without truncation (or repartitioning); with 64 nodes AIRG is 1.8$\times$ faster.

The setup times for $\ell$AIR feature growth in the middle of the hierarchy, similar to that described in \secref{sec:parallel_multigrid}. With AIRG in the previous section, we found better performance from heavily truncated the hierarchy and increasing the number of iterations used in the coarse grid solve. Hence we investigated using 9, 10, 11 and 12 levels with 12 damped Jacobi iterations as a coarse solve with $\ell$AIR. Doing this introduced further growth in the iteration count, giving 18, 20, 21 and 21 iterations and 352, 398, 444 and 444 WUs. This resulted in a reduction of the solve times, despite the higher WUs, confirming the communication bottlenecks on the lower levels. We could not find a coarse grid solver that resulted in both flat iteration counts with this heavily truncated hierarchy and reasonable setup times in \textit{hypre}. This helps show the strength of using a GMRES polynomial as a coarse grid solver in AIRG.

Figures \ref{fig:weak_scaling_solve} \& \ref{fig:weak_scaling_setup} also show the results from Tables \ref{tab:weak_scaling_repartion_truncate} and \ref{tab:weak_scaling_repartion_truncate_heavy} which help demonstrate the impact of repartitioning with ParMETIS onto fewer ranks and truncation with AIRG. We see significantly lower setup and solve times. In particular, with a truncated hierarchy and ParMETIS repartitioning onto fewer ranks, the AIRG solve on 64 nodes is 5.9$\times$ faster than the $\ell$AIR implementation in \textit{hypre} with a heavily truncated hierarchy.
\section{Conclusions}
In this work we examined the parallel scaling of the AIRG reduction multigrid when used to solve Boltzmann transport problems in the streaming limit on unstructured grids. The streaming limit gives a hyperbolic system and hence a semi-coarsening can be produced along the characteristics determined by the angular directions. This gives a coarsening rate of 1/2 (at best) and hence the multigrid hierarchy has many levels. Previously we found distance 2 approximate ideal interpolation gives good results in these problems \cite{Dargaville2024a}, although this long-range interpolation means our coarse grid matrices become increasingly non-local. These factors make the parallelisation of reduction multigrid for hyperbolic systems difficult. 

To achieve good solve and setup times, we found it necessary to repartition our coarse grids onto fewer MPI ranks with ParMETIS and to truncate the multigrid hierarchy. We triggered the repartitioning when the local to non-local nnzs ratio decreased below a fixed threshold; on ARCHER2 we used a value of 2. We applied a GMRES polynomial as the coarse grid solve and found this was very cheap to setup, allowing heavy truncation without affecting the iteration count.

We introduced a two-pass CF splitting, denoted PMISR DDC, which was designed to give a diagonally dominant $\mat{A}_\textrm{ff}$ block while giving good performance in parallel. The first pass uses a maximal independent set algorithm on the symmetrized strong connections, which keeps strong off-diagonal entries out of $\mat{A}_\textrm{ff}$. Many off-diagonal weak connections can still result in rows which are not diagonally dominant, so the second pass then finds the least diagonally dominant rows and sets them as C-points. With our operators we found only a small fraction of rows met this criteria and hence we found a large improvement in the diagonal dominance of $\mat{A}_\textrm{ff}$ with only a small impact on the coarsening rates. We found this PMISR DDC CF splitting produced $\mat{A}_\textrm{ff}$ blocks with the lowest condition numbers, lowest diagonal dominance ratio and lowest worst-case GMRES bounds when compared against Falgout-CLJP, HMIS and a modified PMIS CF splitting algorithm from \textit{hypre}. In parallel this resulted in the fastest solve and setup out of the tested methods. The time spent in the CF splitting for PMISR DDC was between 5--42$\times$ less than any other method. 

When used with our sub-grid scale discretisation, strong scaling results on unstructured grids in 2D showed that using 35k CDOFs/core performed well in the solve with AIRG, though we found reductions in the wall time down to 2.3k CDOFs/core. Using 35k CDOFs/core, we found the grind time for a single V-cycle was as low as 90 ns per DDOF/core.  

Our meshing software limited the weak scaling results we could show and hence we decided to decrease the CDOFs/core and use a low-order angular discretisation (i.e., very few DOFs/spatial node) in order to expose communication bottlenecks. When using 8.8k CDOFs/core, we found that the combination of our PMISR DDC CF splitting, ParMETIS repartitioning onto fewer ranks, truncating our multigrid hierarchy and using a GMRES polynomial as the coarse grid solver resulted in good weak scaling. Moving from 2 to 64 nodes, the iteration count increased from 8 to 9, with the number of work units increasing from 64 to 72. This gives a method that is very close to algorithmically scalable in the hyperbolic limit. The practical scaling in solve time we saw was very close to this ideal, with the time per iteration only increasing by 9.6\%. This gave an overall weak scaling efficiency of 81\% from 2 to 64 nodes on ARCHER2. 

When compared to the $\ell$AIR implementation in \textit{hypre}, we found that AIRG was faster and required less work units to solve. When we introduced ParMETIS repartitioning onto fewer ranks with AIRG and a truncated hierarchy, the solve time was 5.9$\times$ faster than with $\ell$AIR with heavy truncation. A large fraction of this speedup was due to the effect of repartitioning the coarse levels onto fewer MPI ranks, truncation of the hierarchy and the performant coarse grid solver. As such it could be worth introducing these techniques into the \textit{hypre} implementation of $\ell$AIR. Typical elliptic multigrids can use aggressive coarsening to reduce the need for such techniques, but this is not available with reduction multigrids in hyperbolic systems, as the coarsening must be slow. 

As with most algebraic multigrids, we found our setup time grew considerably. We observed a 5.2$\times$ increase in setup time from 2 to 64 nodes, though the overall time was less than the $\ell$AIR implementation in \textit{hypre}. Again repartioning with ParMETIS was vital to ensuring the setup time remained low throughout the hierarchy. We also investigated using the node-aware matrix-vector and matrix-matrix products in the Raptor library. We found if no repartioning was performed the node-aware routines helped reduced the setup time considerably (while also decreasing the solve time), but if repartioning with ParMETIS onto fewer ranks was enabled it was often not worth using the node-aware routines. 

The results from this paper show that it is now possible to achieve good weak scaling and fast solves for hyperbolic Boltzmann transport problems on unstructured grids, without the use of DG discretisations and sweeps. This allows different discretisation and solution methods to be used practically in transport at scale. The weak scaling results were generated with very few DOFs/core, which suggests that good scaling should be possible at higher core counts. Furthermore, the PMISR DDC CF splitting, matrix-vector and matrix-matrix products used are all well suited for use on GPUs. Many aspects of the (expensive) setup can also be re-used in Boltzmann problems, as typically we perform many solves in multi-group and/or eigenvalue problems. Future work will combine the use of our PMSIR DDC CF splitting and AIRG with the additive preconditioning method described in \cite{Dargaville2024a} in order to examine a parallel sweep-free method for solving Boltzmann problems with scattering. 
\section*{Acknowledgements}
This work used the ARCHER2 UK National Supercomputing Service (\url{https://www.archer2.ac.uk}). The authors would like to acknowledge the support of the EPSRC through the funding of the EPSRC grant EP/T000414/1.




\section*{References}
\bibliographystyle{model1-num-names}
\bibliography{bib_library}







\end{document}

%% file: macros.tex
\newcommand{\secref}[1]{Section \ref{#1}}
\newcommand{\fref}[1]{Fig.~\ref{#1}}
\newcommand{\tref}[1]{Table~\ref{#1}}

\renewcommand{\eqref}[1]{Eq.~(\ref{#1})}
\newcommand{\eref}[1]{(\ref{#1})}

\newcommand{\mat}[1]{\textrm{\textbf{#1}}}

\newcommand{\psif}{\psi (\bm{r}, \bm{\Omega}, E)}
\newcommand{\psifd}{\psi (\bm{r}, \bm{\Omega}', E')}

\newcommand{\sn}{S$_\textrm{n}$\,}

\newcommand{\xten}[1]{$\times$ 10$^{\text{#1}}$}

%% file: colours.tex
\definecolor{light}{rgb}{0.8,0.8,0.8}
\definecolor{medium}{rgb}{0.6,0.6,0.6}
\definecolor{dark}{rgb}{0.4,0.4,0.4}
\definecolor{darkmed}{rgb}{0.3,0.3,0.3}
\definecolor{darkest}{rgb}{0.2,0.2,0.2}
\definecolor{Black}{rgb}{0,0,0}
\definecolor{White}{rgb}{1,1,1}
\definecolor{lightpurple}{rgb}{0.78823,0.709803,0.74509}
\definecolor{lightpurpletext}{rgb}{0.788235,0.5529411,0.658823}
\definecolor{skyblue}{rgb}{0.80392,0.866666,0.92941}
\definecolor{skybluetext}{rgb}{0.61568627,0.7647058,0.913725}
\definecolor{darkgreen}{rgb}{0.3137254,0.458823,0.18431}
\definecolor{foliagegreen}{rgb}{0.188,0.415,0.105}
\definecolor{steelbluegrey}{rgb}{0.1961,0.2353,0.2392}
\definecolor{highlightblue}{rgb}{0.4078,0.6431,0.85}
\definecolor{matlabblue}{rgb}{0,0.2705,0.85}
\definecolor{darkred}{rgb}{0.8,0.1725,0}
\definecolor{fireenginered}{rgb}{0.505,0.1411,0}
\definecolor{darkpurple}{rgb}{0.6431,0.3137,0.8509}
\definecolor{gaylordpurple}{rgb}{0.416,0.204,0.549}
\definecolor{deludedorange}{rgb}{0.7409,0.4392,0}
\definecolor{darksalmon}{rgb}{0.9137,0.411,0.706}